\documentclass[11pt,a4paper]{article}

\usepackage{amsmath, amsfonts, amssymb, amsthm}
\usepackage{bbm}
\usepackage{dsfont}
\usepackage{todonotes}
\usepackage{lmodern}
\usepackage[T1]{fontenc}
\usepackage{url}
\usepackage[english]{babel}
\usepackage[numbers,square,sort&compress,longnamesfirst]{natbib}
\usepackage{times,mathptmx}
\usepackage{fullpage}

\newcommand{\zn}[1]{\tr\left(\eta #1\right)}
\numberwithin{equation}{section}
\newcommand{\nm}[1]{\|#1\|}
\newcommand{\A}{\mathcal A}
\newcommand{\1}{\mathds 1}
\newcommand{\R}{\mathbb R}

\newcommand{\N}{\mathbb N}
\newcommand{\Z}{\mathbb Z}

\newcommand{\B}{\mathcal B}
\newcommand{\E}{\mathbb E}
\newcommand{\s}{\mathbb S}
\renewcommand{\vec}{\mathrm{vec}}
\DeclareMathOperator{\tr}{tr}
\renewcommand{\Re}{\re}
\newcommand{\p}{\mathbb P}
\newcommand{\q}{\mathbb Q}

\newcommand{\abs}[1]{\,\left| \!\, #1 \!\, \right|\,}

\newcommand{\sdd}{{\mathbb S_d^+}}

\newcommand{\cb}{\mathcal C_b}
\newcommand{\co}{\mathcal C_0}
\newcommand{\re}{\operatorname{Re}}

\newtheorem{theorem}{Theorem} \numberwithin{theorem}{section} 
\newtheorem{definition}[theorem]{Definition}
\newtheorem{bem}[theorem]{Remark}\newtheorem{remark}[theorem]{Remark}\newtheorem{lemma}[theorem]{Lemma}
\newtheorem{prop}[theorem]{Proposition}
\newtheorem{cor}[theorem]{Corollary}
\newtheorem{example}[theorem]{Example}

\title{Geometric Ergodicity of the multivariate COGARCH(1,1) Process }
\author{Robert Stelzer\thanks{Institute of Mathematical Finance, Ulm University, Helmholzstra\ss e 18, 
D-89069 Ulm, E-mail: robert.stelzer@uni-ulm.de} \and Johanna Vestweber\thanks{Institute of Mathematical Finance, Ulm University, Helmholzstra\ss e 18, 
D-89069 Ulm, E-mail: johanna.vestweber@uni-ulm.de}}
\date{}

\begin{document}
\maketitle

\begin{abstract}
For the multivariate COGARCH(1,1) volatility process we show sufficient conditions for the existence of a unique stationary distribution, for the geometric ergodicity and for the finiteness of moments of the stationary distribution by a Foster-Lyapunov drift condition approach. The test functions used are naturally related to the geometry of the cone of positive semi-definite matrices and the drift condition is shown to be satisfied if the drift term of the defining stochastic differential equation is sufficiently ``negative''.  We show easily applicable sufficient conditions for the needed irreducibility and aperiodicity of the volatility process living in the cone of positive semidefinite matrices, if the driving L\'evy process is a compound Poisson process. 
\end{abstract}

\begin{tabbing}
\end{tabbing}

\vspace{0.3cm}


\section{Introduction}
General autoregressive conditionally heteroscedastic (GARCH) time series models, as introduced in \cite{bollerslev1986}, are of high interest for financial economics. They capture many typical features of observed financial data, the so-called stylized facts (see \cite{Guillaumeetal1997}). A continuous time extension, which captures the same stylized facts as the discrete time GARCH model, but can also be used for irregularly-spaced and high-frequency data, is the COGARCH process, see e.g. \cite{KLM2004, KLM2006, BCL2006}. The use in financial modelling is studied e.g. in \cite{KlueppelbergSzimayerMaller2011,Panayotov2005,BayraciUenal2014} and the statistical estimation in \cite{Haugetal2007,BibbonaNegr2015,Malleretal2006}, for example. Furthermore, an asymmetric variant is proposed in \cite{BehmeKlueppelbergMayr2014} and an extension allowing for more flexibility in the autocovariance function in \cite{behmeKlueppelbergChong2015}.  

To model and understand the behavior of several interrelated time series as well as to price derivatives on several underlyings or to assess the risk of a portfolio multivariate models for financial markets are needed. The fluctuations of the volatilities and correlations over time call for employing stochastic volatility models which in a multivariate set-up means that one has to specify a latent process for the instantaneous covariance matrix. Thus, one needs to consider appropriate stochastic processes in the cone of positive semi-definite matrices.  Many popular multivariate stochastic volatility models in continuous time, which in many financial applications is preferable to modelling in discrete time, are of an affine type, thus falling into the framework of \cite{CFMT,MayerhoferStelzerVestweber2018}.  Popular examples include the Wishart (see e.g. \cite{BaldeauxPlaten2013,bru,GnoattoGrasselli2014}) and the Ornstein-Uhlenbeck type stochastic volatility model (see \cite{MuhleKarbePfaffelStelzer2009,BenthVos2013,GranelliVeraart2016}, for example). 
 Thus the price processes have two driving sources of randomness and the tail-behavior of their volatility process is typically equivalent to the one of the driving noise (see \cite{Moseretal2010,Fasenetal2004}). A very nice feature of GARCH models is that they have only one source of randomness and their structure ensures heavily-tailed stationary behavior even for very light tailed driving noises (\cite{FernandezMuriel2009,Basraketal2002}).
In discrete time one of the most general multivariate GARCH versions (see \cite{Bauwensetal2006,FrancqZakoian2010} for an overview) is the BEKK model, defined in \cite{englekroner}, and the multivariate COGARCH(1,1) (shortly MUCOGARCH(1,1)) process introduced and studied in \cite{S2010}, is the continuous time analogue, which we are investigating further in this paper.

The existence and uniqueness of a stationary solution as well as the convergence to the stationary solution is of high interest and importance. Geometric ergodicity ensures fast convergence to the stationary regime in simulations and paves the way for statistical inference. By the same argument as in \cite[Proof of Theorem 4.3, Step 2]{masuda2004} geometric ergodicity and the existence of some $p$-moments of the stationary distribution provide exponential $\beta$-mixing for Markov processes. This in turn can be used to show a central limit theorem for the process, see for instance \cite{doukhan}, and so allows to prove for example asymptotic normality of estimators (see e.g. \cite{Haugetal2007,BibbonaNegr2015} in the context of univariate COGARCH(1,1) processes). In a similar way \cite{doRegoSousaKlueppelbergStelzer2018,doRegoSousa2019} employ the results of the present paper to analyse moment based estimators of the parameters of MUCOGARCH(1,1) processes.

In many applications involving time series (multivariate) ARMA-GARCH models (see e.g. \cite{FrancqetalARMAGARCH2004,LingetMcAleer2003,ComteLieberman2003}) turn out to be adequate and geometric ergodicity is again key to understand the asymptotic behaviour of statistical estimators. In continuous time a promising analogue currently investigated in \cite{RezapourStelzer2017} seems to be a (multivariate)  CARMA process (see e.g. \cite{Brockwell2009,Fasen:Kimmig:2016,marquardt:stelzer:2007}) driven by a (multivariate) COGARCH process. The present paper also lays  foundations for the analysis of such models.

For the univariate COGARCH process geometric ergodicity was shown by \cite{fasen2010} and \cite{fuchs} discussed it for the BEKK GARCH process. In \cite{S2010} for the MUCOGARCH process sufficient  conditions for the  existence of a stationary distribution are shown by tightness arguments, but the paper failed to establish uniqueness or convergence to the stationary distribution. In this paper we deduce under the assumption of irreducibility sufficient conditions for the uniqueness of the stationary distribution, the convergence to it with an exponential rate and some finite $p$-moment of the stationary distribution of the MUCOGARCH volatility process $Y$. To show this we use the theory of Markov process, see e.g. \cite{MT2, DMT1995}. A further result of this theory is, that our volatility process is positive Harris recurrent. If the driving L\'evy process is a compound Poisson process, we show easily applicable conditions ensuring irreducibility of the volatility process in the cone of positive semidefinite matrices. 

Like in the discrete time BEKK case the non-linear structure of the SDE prohibits us from using well-established results for random recurrence equations like in the one-dimensional case and due to the rank one jumps establishing irreducibility is a very tricky issue. To obtain the latter \cite{fuchs} in discrete time   used techniques from algebraic geometry (see also \cite{Boussama1998th2}) whereas we use a direct probabilistic approach playing the question back to the  existence of a density for a Wishart distribution. However, we restrict ourselves to processes of order (1,1) while in the discrete time BEKK case general orders were considered. The reason is that on the one hand order (1,1) GARCH processes seem sufficient in most applications and on the other hand multivariate COGARCH(p,q) processes can be defined in principle (\cite[Section 6.6]{S2007}), but no reasonable conditions on the possible parameters are known. Already in the univariate case these conditions are quite involved (cf. \cite[Section 5]{BCL2006}, \cite{TsaiChan2009}). On the other hand we look at the finiteness of an arbitrary $p$-th moment (of the volatility process) and use drift conditions related to it, whereas \cite{fuchs} only looked at the first moment for the BEKK case. In contrast to \cite{S2010} we avoid any vectorizations, work directly in the cone of positive semi-definite matrices and use test functions naturally in line with the geometry of the cone.

After a brief summary of some preliminaries, notations and L\'evy processes, the remainder of the paper is organized as follows: In Section 3 we recall the definition of the MUCOGARCH(1,1) process and some of its properties of relevance later on. In Section 4 we present our first main result: sufficient conditions ensuring the geometric ergodicity of the volatility process $Y$. Furthermore, we compare the conditions for geometric ergodicity to previously known conditions for (first order) stationarity.  Moreover, we discuss the applications of the obtained results and illustrate them by exemplary simulations. In Section \ref{section:irredresults} we establish sufficient conditions for the  irreducibility and aperiodicity of $Y$ needed to apply the previous results on geometric ergodicity. Section \ref{section:proofsMarkov} first gives a brief repetition of the Markov theory we use and the proofs of our results are developed.

\section{Preliminaries} 
Throughout we assume that all random variables and processes are defined on a given filtered probability space 
$(\Omega,\mathcal{F}, \p, (\mathcal{F}_t)_{t\in \mathcal{T}})$ with $\mathcal{T}=\N$ in the discrete time case and $\mathcal{T}=\R^+$ in the 
continuous one. Moreover, in the continuous time setting we assume the usual conditions (complete, right continuous filtration) to be satisfied. 

For Markov processes in discrete 
and continuous time we refer to \cite{MT} and respectively \cite{dynkinI, applebaum}.
A summary of the most relevant notions and results from Markov processes is given in Section \ref{chaptermarkov} for the convenience of the reader.

\subsection{Notation}

The set of real $m \times n$ matrices is denoted by $M_{m,n}(\R)$ or only by $M_n(\R)$ if $m=n$. For the invertible $n \times n$ matrices
we write $GL_n(\R)$. The linear subspace of symmetric matrices we denote by $\s_n$, by $\s_n^+$ the closed cone of   positive semi-definite 
matrices and the open cone of positive definite matrices by $\s_n^{++}$. Further we denote by $I_n$ the $n \times n$ identity matrix. \\
We introduce the natural ordering on $\s_n$ and denote it by $\preceq$, that is for $A,B  \in \s_n$ it holds $A\preceq B \ \Leftrightarrow \  B-A \in \s_n^+$.
The tensor (Kronecker) product of two matrices $A,B$ is written as $A\otimes B$. $\vec$ denotes the well-known vectorization operator that maps 
the $n\times n$ matrices to $\R^{n^2}$ by stacking the columns of the matrices below another.  The spectrum of a matrix is denoted by $\sigma(\cdot)$ and the spectral radius by $\rho(\cdot)$. For a matrix with only real eigenvalues $\lambda_{max}(\cdot)$ and $\lambda_{min}(\cdot)$ denote the largest and the smallest eigenvalue. $\re(x)$ is the real part of a complex number. Finally, $A^\top$ is the transpose of a matrix $A\in M_{m,n}(\R)$.

By $\| . \|_\text{2}$ we denote both the Euclidean norm for vectors and the 
corresponding operator norm for matrices and by $\|.\|_\text{F}$ the 
Frobenius norm for matrices.

Furthermore, we employ an intuitive notation with respect to the (stochastic) integration with matrix-valued integrators referring to any of 
the standard texts (e.g. \cite{Protter2010}) for a comprehensive treatment of the theory of stochastic integration. For an $M_{m,n}(\R)$-valued L\'evy process $L$, and $M_{d,m}(\R)$ resp.\ $M_{n,p}(\R)$- valued  processes $X,Y$ integrable with respect to $L$, the term $\int_0^t X_s \, dL_s Y_s$ is to be understood as the $d\times p$ (random) matrix with $(i,j)$-th entry  $\sum_{k=1}^m \sum_{l=1}^n \int_0^t X_s^{ik} \, dL_s^{kl} Y_s^{lj}$.
If $(X_t)_{t\in\R^+}$ is a semi-martingale in $\R^m$ and $(Y_t)_{t\in\R^+}$ one in $\R^n$ then the quadratic variation $([X,Y]_t)_{t\in\R^+}$ is 
defined as the finite variation process in $M_{m,n}(\R)$ with components $[X,Y]_{ij,t}=[X_i,Y_j]_t$ for $t\in\R^+$ and $i=1,\ldots,m$, $j=1,\ldots,n$. 


\subsection{L\'evy processes}
Later on we use L\'evy processes (see e.g. \cite{applebaum,sato}) in $\R ^d$ and in the symmetric matrices $\s_d$. 

We consider a L\'evy process $L=(L_t)_{t\in\R ^+}$ (where $L_0=0$ a.s.)
in $\R ^d$  determined by its characteristic function in the 
L\'evy-Khintchine form
$E\left[e^{i\langle u,{L}_t\rangle}\right]=\exp\{t\psi_L(u)\}$ for $t\in\R ^+$ with
\begin{equation*} 
\psi_L(u)=i\langle\gamma_L,u\rangle-\frac{1}{2}\langle u,\tau_Lu\rangle+\int\limits_{\R ^d}
(e^{i\langle u,x\rangle}-1-i\langle u,x\rangle I_{[0,1]}(\{\nm{x}))\,\nu_L(dx),\quad u\in\R ^d,
\end{equation*}
where $\gamma_L\in\R ^d$, $\tau_L\in\s_d^+$ and the L\'evy measure $\nu_L$ is a measure on $\R ^d$ satisfying
$\nu_L(\{0\})=0$ and  $\int_{\R ^d}(\nm{x}^2\wedge 1)\,
\nu_L(dx)<\infty.$ If $\nu_L$ is a finite measure, $L$ is a compound Poisson process.
Moreover, $\langle\cdot,\cdot\rangle$ denotes the usual Euclidean scalar product on $\R ^d$.

We always assume $L$ to be c\`adl\`ag and denote its jump measure by $\mu_L$, i.e. $\mu_L$ is the Poisson random measure on $\R ^+\times \R ^d\setminus \{0\}$ given by $ \nonumber
\mu_L(B)=\sharp\{s\ge 0:\;(s,L_s-L_{s-})\in B\}
$ for any
measurable set $B\subset \R ^+\times\R ^d\setminus \{0\}.$ Likewise, $\tilde \mu_L(ds,dx)=\mu_L(ds,dx)-ds\nu_L(dx) $ denotes the compensated jump measure.

Regarding matrix-valued L\'evy processes, we will only encounter matrix subordinators (see \cite{BarndorffetPerez2005}), i.e. L\'evy processes with paths in $\s_d^+$. Since matrix subordinators are  of finite variation and $\tr(X^*Y)$ (with $X,Y\in
\s_d$ and $\tr$ denoting the usual trace functional) defines a
scalar product on $\s_d$  linked to the Euclidean scalar product on $\R ^{d^2}$ via  $\tr(X^*Y)=\vec(X)^*\vec(Y)=\langle \vec(Y), \vec(X)\rangle$, the characteristic function of a matrix subordinator can be represented as
\begin{align*}
E\left(e^{i\tr(L_t^*Z)}\right)&=\exp\left(t\psi_L(Z)\right),\, Z\in \s_d,\,\mbox{where}\,
\psi_L(Z)=i \tr(\gamma_L    Z)+\int_{\s_d^+}(e^{i\tr(XZ)}-1)\nu_L(dX)
\end{align*}
with drift $\gamma_L\in\s_d^+$ and L\'evy measure $\nu_L$.

The discontinuous part of the quadratic variation of any L\'evy process $L$ in $\R ^d$
\[
[L,L]^d_t:=\int_0^t\int_{\R ^d}xx^*\mu_L(ds,dx)=\sum_{0\leq s\leq t}(\Delta L_s)(\Delta L_s)^*,
\]
is  a matrix subordinator with drift zero and L\'evy measure 
$
\nu_{[L,L]^d}(B)=\int_{\R ^d}I_B(xx^*)\nu_L(dx)
$
for all Borel sets $B\subseteq \s_d$.
\section{Multivariate COGARCH(1,1) process}\label{sectiondefmucogarch}
In this section we present the definition of the MUCOGARCH(1,1) process and some relevant  properties mainly based on \cite{S2010}.

\begin{definition}[MUCOGARCH(1,1), {\cite[Definition 3.1]{S2010}},] Let $L$ be an $\R^{d}$-valued L\'evy process, $A,B \in M_{d}(\R)$ and $C \in \s_{d}^{++}$. 
The \textbf{MUCOGARCH(1,1) process} $G=(G_{t})_{t \ge 0}$ is defined as the solution of 
\begin{align}
 dG_{t} &= V_{t-}^{\frac 1 2}dL_{t} \\
 \label{defv}V_{t} &= Y_{t} + C \\
 \label{defy} dY_{t} &= (BY_{t-} + Y_{t-}B^{\top}) dt + A V_{t-}^{\frac 1 2} d [L,L]^{d}_{t}V_{t-}^{\frac 1 2} A^{\top},
\end{align}
with initial values $G_{0} \in \R^{d}$ and $Y_{0} \in \s_{d}^{+}$.  

The process $Y=(Y_{t})_{t\ge 0}$ is called \textbf{MUCOGARCH(1,1) volatility process}. 
\end{definition}
Since we only consider MUCOGARCH(1,1) processes, we 
often simply write MUCOGARCH. 

Equations \eqref{defv} and \eqref{defy} directly give us an SDE for the covariance matrix process $V$:
\begin{equation}
 dV_t = (B(V_{t-} -C)+(V_{t-}-C)B^{\top})dt + AV_{t-}^{\frac 1 2} d[L,L]^d_t 
V_{t-}^{\frac 1 2} A^{\top}.\label{sdev}
\end{equation}

Provided $\sigma(B) \subset (-\infty,0) + i\R$, we see that $V$, as long as no jumps occur,  
returns to the level $C$ at an exponential rate determined by $B$. Since all 
jumps are positive semidefinite, $C$ is not a mean level but a lower bound.

To have the MUCOGARCH process well-defined, we have to know that a unique 
solution of the SDE system exists and the solution of $Y$ (and $V$) does 
not leave the set $\s_d^+$. In the following we always understand that our processes live on $\s_d$. Since $\s_d^{++}$ is an open subset of $\s_d$, we now are in the most natural setting for SDEs and we get:

\begin{prop}[{\cite{S2010}}, Theorems 3.2, 3.6] \label{theoremvolatility}
Let $A,B\in M_d(\R)$, $C\in \s_d^{++}$  and $L$ be a $d$-dimensional  L\'evy process. 
Then the SDE (\ref{defy}) with  initial value $Y_0\in\s_d^+$ has a unique positive semi-definite solution $(Y_t)_{t\in\R^+}$.
The solution $(Y_t)_{t\in\R^+}$ is locally bounded and of finite variation. 

Moreover, it satisfies \begin{equation}
Y_t=e^{Bt}Y_0e^{B^\top t}+\int_0^te^{B(t-s)}A(C+Y_{s-})^{1/2}d[L,L]_s^d (C+Y_{s-})^{1/2}A^\top e^{B^\top(t-s)}\label{volterra}
\end{equation}
for all $t \in \R^+$ and thus $Y_t\succeq e^{Bt}Y_0e^{B^\top t}$ 
for all $t\in\R^+$.
\end{prop} 
In particular, whenever \eqref{defy} is started with an initial value $Y_0\in\s_d^+$ (or \eqref{sdev} with $V_0\succeq C$) the solution stays in $\s_d^+$ ($\s_d^++C$) at all times. This can be straightforwardly seen from \eqref{volterra} and the fact that for any $M\in M_d(\R)$ maps of the form $X\mapsto MXM^\top$ map $\s_d^+$ into itself (see \cite{S2010} for a more detailed discussion).

At first sight, it appears superfluous to introduce the process $Y$ instead of directly working with the process $V$. However, in the following it is more convenient to work with $Y$, as then the matrix $C$ does not appear in the drift and the state space of the Markov process analysed is the cone of positive semi-definite matrices itself and not a translation of this cone.

\begin{prop}[Markov properties]\label{th:markov}
\begin{enumerate}
\item[(i)] The MUCOGARCH process $(G,Y)$ as well as its volatility process $Y$ alone are  
temporally homogeneous  strong Markov processes on $\R^d\times\s_d^{+}$ and 
$\s_d^{+}$, respectively, and they have the weak $\cb$-Feller property.
\item[(ii)] $Y$ is non-explosive and has the weak $\co$-Feller property. Thus it is a Borel right process.
\end{enumerate}
\end{prop}
\begin{proof}
	\cite[Theorem 4.4]{S2010} is (i).
	
	 $\co$-Feller property of $Y$: Let $f\in \co(\s_d^+)$. We have to show that $\forall 
	t\ge0$
	$ 
	\p^t f(x) \to 0,\  \   \text{ for } x\to\infty,
	$
	where we understand $x \to \infty$ in the sense of $\|x\|_2\to \infty$.
	Since $\p^t f(x)=\E(f(Y_t(x)))$, where $x$ denotes the starting point $Y_0=x$, it is enough to show that $Y_t$ goes to infinity for $x \to \infty$:
	\begin{align*}
	\|Y_t\|_2 \ge & \| e^{Bt}Y_0 e^{B^\top t} \|_2 	\geq \| e^{-Bt} \|^{-2}_2 \|x\|_2 \to \infty \text{ for } \|x\|_2 \to \infty. 
	\end{align*}
	
As argued in Section	\ref{chaptermarkov} any  $\co$-Feller process is a Borel right process and the non-explosivity property is shown in the proof of Theorem 6.3.7 in 
	\cite{S2007}. 
\end{proof}
\begin{bem}
For $\R^d$-valued solutions to L\'evy-driven stochastic differential equations 
$dX_t = \sigma(X_{t-})dL_t $
\cite{kuhn2016} gives necessary and sufficient conditions for the rich Feller property, which includes the $\co$-Feller property, if $\sigma$ is continuous and of (sub-)linear growth. But since our direct proof is quite short, we prefer it instead of trying to adapt the result of \cite{kuhn2016} to our state space. 
\end{bem}

\section{Geometric ergodicity of the MUCOGARCH volatility process $Y$} 
\label{section:mainresults}

In Theorem \cite[Theorem 4.5]{S2010} sufficient conditions for the existence of a stationary distribution for the volatility process $Y$  or $V$ (with certain moments finite) are shown, but neither the uniqueness of the stationary distribution nor that it is a limiting distribution are obtained. Our main theorem now gives sufficient conditions for geometric ergodicity and thereby for the existence of a unique stationary distribution to which the transition probabilities converge exponentially fast in total variation (and in stronger norms).

On the proper closed convex cone $\sdd$ the trace $\tr:  \sdd \to \R^+$ is well-known to define a norm which is also a linear functional. So do actually all the maps $\sdd \to \R^+,\, X\mapsto \tr(\eta X)$ with $\eta\in\s_d^{++}$, as $\s_d^+$ is also generating, self-dual and has interior $\s_d^{++}$ (cf. \cite{faraut1994analysis}, for instance). The latter can also be easily seen using the following Lemma which is a consequence of \cite{Richter1958}.
\begin{lemma}\label{lem:eigtrace}
Let $X\in \s_d$, $Y\in \s_d^+$. Then
\[
\lambda_{min}(X)\tr(Y)\leq \tr(XY)\leq \lambda_{max}(X)\tr(Y).
\]	
\end{lemma}

To be in line with the geometry of $\s_d^+$ we thus use the above norms to define appropriate test functions and look at the trace norm for the finiteness of moments (which, of course, is independent of the actually employed norm).

For $p=1$ and $p \ge2$ it is shown in \cite[Proposition 4.7]{S2010}, that the finiteness of $\E(\|L_1\|_2^{2p})$ and $\E(\|Y_0\|_2^{p})$ implies the finiteness of $\E(\|Y_t\|_2^{p})$ for all $t$. We improve this to all $p >0$. 
\begin{lemma}[Finiteness of moments] \label{lemma:moments}
	Let $Y$ be a MUCOGARCH volatility proces and  $p >0$. If $\E(\tr(Y_0)^{p}) < \infty$, $\int_{\|y\|_2\leq 1} \| y \|_2 ^{2\wedge2p} \nu_L(dy) < \infty$  and $\E(\|L_1\|_2^{2p}) < \infty$, then $\E(\tr(\eta Y_t)^{p})< \infty$ for all $t \ge 0$,  $\eta\in\s_d^{++}$ and $t \mapsto \E(\tr(\eta Y_t)^{p})$ is locally bounded. 
\end{lemma}
Observe that $\E(\|L_1\|_2^{2p}) < \infty$ is equivalent to $\int_{\|y\|_2\geq 1} \| y \|_2 ^{2p} \nu_L(dy) < \infty$ and that  $\int_{\|y\|_2\leq 1} \| y \|_2 ^{2\wedge2p} \nu_L(dy) < \infty$ is always true for $p\geq 1$ and otherwise means that the $2p$-variation of $L$ has to be finite.
\begin{theorem}[Geometric ergodicity] \label{theoremge}

Let $Y$ be a MUCOGARCH volatility process  which is $\mu$-irre\-ducible with the support of $\mu$ having 
non-empty interior and aperiodic.

If one of the following conditions is satisfied
\begin{enumerate}
	\item[(i)] setting $p=1$ there exists an $\eta\in \s_d^{++}$ such that $\int_{\|y\|_2\geq 1}\|y\|_2^{2}\nu_L(dy)<\infty$ and 
	\begin{equation} \label{eq:cond1}
	\eta B+B^\top \eta+ A^\top\eta A\left\|\int_{\R^{d}}  yy^\top  \nu_L(dy)\right\|_2\in-\s_d^{++},
	\end{equation}
	\item[(ii)] setting $p=1$ there exists an $\eta\in \s_d^{++}$ such that $\int_{\|y\|_2\geq 1}\|y\|_2^{2}\nu_L(dy)<\infty$ and 
	\begin{equation} \label{eq:cond2}
	\eta B+B^\top \eta+ \lambda_{max}(A^\top\eta A)\int_{\R^{d}} yy^\top\nu_L(dy)\in-\s_d^{++},
	\end{equation}
	\item[(iii)] there exist a $p\in (0,1]$  and an $\eta\in \s_d^{++}$ such that 	$\int_{\R^d} \| y \|_2 ^{2p} \nu_L(dy) < \infty$ and 
	\begin{equation} \label{eq:assumptionp<1}\int_{\R^{d}} \left( \left( 1+K_{\eta,A}\|y\|^2_2\right)^p - 1\right)\nu_L(dy)
	+K_{\eta,B} p  < 0,
	\end{equation}
	where $K_{\eta,B}:=\max\limits_{x\in\s_d^+,\tr(x)=1 }\frac{\tr((\eta B+B^\top\eta)x)}{\tr(\eta x)}$ and $K_{\eta,A}:=\max\limits_{x\in\s_d^+,\tr(x)=1 }\frac{\tr(A^\top\eta Ax)}{\tr(\eta x)}$,
 \item[(iv)] there exist a $p\in [1,\infty)$ and an $\eta\in \s_d^{++}$ such that $\int_{\|y\|_2\geq 1}\|y\|_2^{2p}\nu_L(dy)<\infty$ and 
	\begin{equation} \label{eq:assumptionpge1}\int_{\R^{d}} \left(2^{p-1} \left( 1+K_{\eta,A}\|y\|^2_2\right)^p - 1\right)\nu_L(dy)
	+K_{\eta,B} p  < 0,
\end{equation}	
where $K_{\eta,B}, K_{\eta,A}$ are as in (iii),
 \item[(v)] there exists a $p\in [1,\infty) $ and an $\eta\in \s_d^{++}$ such that $\int_{\|y\|_2\geq 1}\|y\|_2^{2p}\nu_L(dy)<\infty$ and 
 \begin{equation} \label{eq:assumptionv}
 \max\{2^{p-2},1\}K_{\eta,A}\int_{\R^{d}}\|y\|^2_2\left(  1 +\|y\|^2_2K_{\eta,A}\right)^{p-1}\nu_L(dy)+K_{\eta,B}    < 0
 \end{equation}
 where $K_{\eta,B}, K_{\eta,A}$ are as in (iii),
 \end{enumerate}
then a  unique stationary distribution for the MUCOGARCH(1,1)
volatility process $Y$ exists, $Y$ is positive Harris recurrent, geometrically ergodic (even $\tr(\eta\cdot)^p+1$ uniformly ergodic) and the 
stationary distribution has a finite $p$-th moment.
\end{theorem}
\begin{remark}\label{rem44}
	\begin{itemize}
	\item[(i)] For $p=1$ the cases (i) and (ii) give us quite strong conditions comparable to the ones known for affine processes in $\sdd$ (cf. \cite{MayerhoferStelzerVestweber2018}). For $p\not = 1$ it seems that the non-linearity of our SDE implies that we need to use inequalities that are somewhat crude.
	\item[(ii)] Condition \eqref{eq:assumptionpge1} demands that the driving L\'evy process is compound Poisson for $p>1$. To overcome this restriction is the main motivation for considering case (v).
	\item[(iii)] For $p=1$ the Conditions  \eqref{eq:assumptionp<1}, \eqref{eq:assumptionpge1}, \eqref{eq:assumptionv} agree. 
	
	In dimension one they also agree with Conditions \eqref{eq:cond1} and \eqref{eq:cond2}. Moreover, then the above sufficient conditions agree with the necessary and sufficient conditions of \cite[Lemma 4.1]{KLM2004} for a univariate COGARCH(1,1) process to have a stationary distribution with finite first moments, as follows from \cite[p. 88]{S2010}.
	
	\item[(iv)] 	 Observing that \[
	K_{\eta,A}\int_{\R^{d}}\|y\|^2_2\nu_L(dy)+K_{\eta,B} \geq \max\limits_{x\in\s_d^+,\tr(x)=1 }\frac{\tr((\eta B+B^\top\eta+A^\top\eta A\left\|\int_{\R^{d}}  yy^\top  \nu_L(dy)\right\|_2)x)}{\tr(\eta x)},
	\]
	the self-duality of $\sdd$ shows that for $p=1$ the equivalent Conditions \eqref{eq:assumptionp<1}, \eqref{eq:assumptionpge1}, \eqref{eq:assumptionv} imply  that \eqref{eq:cond1} holds. 
	
	Conversely the upcoming Examples \ref{condexamp1}, \ref{condexamp2} show that \eqref{eq:cond1} and  \eqref{eq:cond2} are less restrictive than the equivalent Conditions \eqref{eq:assumptionp<1}, \eqref{eq:assumptionpge1}, \eqref{eq:assumptionv} and that \eqref{eq:cond1} does not imply  \eqref{eq:cond2} and vice versa. 
\item[(v)] 	Arguing as in \cite[Lemma 4.1]{KLM2004}, one sees that if \eqref{eq:assumptionp<1} is satisfied for some $p>0$ it is also satisfied for all  smaller ones. Using similar elementary arguments, the same can be shown for \eqref{eq:assumptionpge1} and for \eqref{eq:assumptionv} this property is obvious.
	
	Note, however, that $\int_{\R^d} \| y \|_2 ^{2p} \nu_L(dy) < \infty$ for some $0<p\leq 1$ does not imply that $\int_{\R^d} \| y \|_2 ^{2\bar p} \nu_L(dy)$ is finite for $0<\bar p \leq p$.
	\item[(vi)] From the exercise on p. 98 of \cite{Hornetal1991} (taking e.g. $A=\begin{pmatrix} 1 & 0 \\ 10 & 1 \end{pmatrix}$ there) we see immediately that if \eqref{eq:cond1},  \eqref{eq:cond2}, \eqref{eq:assumptionp<1}, \eqref{eq:assumptionpge1}, or \eqref{eq:assumptionv} is satisfied for one $\eta\in\s_d^{++}$ it can be violated for $\bar\eta\in\s_d^{++}, \bar\eta\neq \eta.$ Hence, the possibility to choose $\eta\in\s_d^{++}$ freely is important.
	\item[(vii)] In the Case (iii) $\int_{\R^{d}} K_{\eta,A}^p\|y\|^{2p}_2\nu_L(dy)
	+K_{\eta,B} p  < 0,
$ implies \eqref{eq:assumptionp<1}, as $(x+y)^p-x^p\leq y^p$ for $x,y\geq 0$ and $p\in(0,1]$.
		\end{itemize}
\end{remark}
The above Conditions \eqref{eq:assumptionp<1}, \eqref{eq:assumptionpge1}, \eqref{eq:assumptionv} appear enigmatic at a first sight. However, essentially they are related to the drift being negative in an appropriate way.
\begin{lemma}\label{lem:etaB}
	\begin{itemize}
		\item[(i)] If one of the Conditions \eqref{eq:cond1}, \eqref{eq:cond2}, \eqref{eq:assumptionp<1}, \eqref{eq:assumptionpge1}, \eqref{eq:assumptionv} is satisfied, then $\eta B+B^\top\eta\in-\s_d^{++}$.
		\item[(ii)] If one of the Conditions \eqref{eq:cond1}, \eqref{eq:cond2}, \eqref{eq:assumptionp<1}, \eqref{eq:assumptionpge1}, \eqref{eq:assumptionv} is satisfied, then $\Re(\sigma(B))<0$.
		\item[(iii)] If $\Re(\sigma(B))<0$, then there exists an $\eta\in \s_d^{++}$ such that $\eta B+B^\top\eta\in-\s_d^{++}$.
	\end{itemize}
\end{lemma}
\begin{proof}
(i)	This is obvious for \eqref{eq:cond1}, \eqref{eq:cond2}.
	In the other cases we get that $K_{\eta,B}<0$ must hold. Therefore $\max\limits_{x\in\s_d^+,\tr(x)=1 }\tr((\eta B+B^\top\eta)x)<0$. By the self-duality of $\sdd$ this shows $\eta B+B^\top\eta\in-\s_d^{++}$.

(ii) and (iii) now follow  from \cite[Theorem 2.7, Corollary 5.1]{MayerhoferStelzerVestweber2018}.
\end{proof}
So in the end what we demand is that the drift is ``negative'' enough to compensate for a positive effect from the jumps.  The Conditions \eqref{eq:cond1}, \eqref{eq:cond2} can be related to eigenvalues of linear maps on $\s_d$.
\begin{lemma}\label{lem:linmaps}
	\begin{itemize}
		\item[(i)] Condition \eqref{eq:cond1} holds if and only if the linear map $\s_d\to \s_d$, $X\mapsto XB+B^\top X +A^\top X A \left\|\int_{\R^{d}} yy^\top \nu_L(dy)\right\|_2 $ has only eigenvalues with strictly negative real part.
		\item[(ii)] Condition \eqref{eq:cond2} holds if the linear map $\s_d\to \s_d$, $X\mapsto XB+B^\top X +\tr(A^\top X A)\int_{\R^{d}}  yy^\top  \nu_L(dy)$ has only eigenvalues with strictly negative real part.
	\end{itemize}
\end{lemma}	
\begin{proof}
	(i)	Follows from \cite[Theorem 2.7]{MayerhoferStelzerVestweber2018}, as it is easy to see that the map is quasi monotone increasing and as a linear map and its adjoint have the same eigenvalues. 
	
	(ii) Follows analogously after noting that $0\leq \lambda_{max}( A^\top X A)\leq \tr(A^\top X A)$.
\end{proof}

\begin{bem}[Relation to first order stationarity]
	We say that a process is first order stationary, if it has finite first moments at all times, if the first moment converges to a finite limit independent of the initial value as time goes to infinity and if the first moment is constant over time when the process is started at time zero with an initial value whose  first moment equals the limiting value.
	
	According to \cite[Theorem 4.20, its proof and Remark 4.9]{S2010} sufficient conditions for asymptotic first-order stationarity of the MUCOGARCH volatility $Y$ are:
	\begin{enumerate}
		\item[(a)] there exists a constant $\sigma_L \in \R^+$ such that $\int_{\R^d} y y^\top \nu_L(dy) = \sigma_L I_d$,
		\item[(b)] $\sigma(B) \subset (-\infty, 0) + i\R$,
		\item[(c)] $\sigma \left( B \otimes I + I \otimes B + \sigma_L(A \otimes A)  \right) \subset (-\infty, 0) + i\R. $ 
	\end{enumerate}
	An inspection of the arguments given there shows that (c) only needs to hold for the linear operator $B \otimes I + I \otimes B + \sigma_L(A \otimes A)$  restricted to the set $\vec(\s_d)$. Under (a) $\left\|\int_{\R^{d}} yy^\top \nu_L(dy)\right\|_2=\sigma_L$ and devectorizing thus shows that $B \otimes I + I \otimes B + \sigma_L(A \otimes A)$ is just the linear operator in Lemma \ref{lem:linmaps} (i). 
	
	Hence, under (a) our Condition \eqref{eq:cond1} implies (b) and (c). So our conditions for geometric ergodicity with a finite first moment  are certainly not worse than the previously known conditions for just first order stationarity. 
\end{bem}

The constants $K_{\eta,B},\,K_{\eta,A}$ are related to changing norms and may be somewhat tedious to obtain. The following lemma follows immediately from Lemma \ref{lem:eigtrace} and implies that we can replace them by eigenvalues in the Conditions  \eqref{eq:assumptionp<1}, \eqref{eq:assumptionpge1}, \eqref{eq:assumptionv}.
\begin{lemma}
It holds that:
\begin{enumerate}
	\item $K_{\eta,B}\leq\frac{\lambda_{max}(\eta B+B^\top\eta)}{\lambda_{min}(\eta)}$.
	\item $K_{\eta,A}\leq\frac{\lambda_{max}(A^\top\eta A)}{\lambda_{min}(\eta)}$.
	\item $K_{I_d,B}={\lambda_{max}(B+B^\top)}$.
	\item $K_{I_d,A}={\lambda_{max}(A^\top A)}$.
\end{enumerate}
\end{lemma}	
\begin{bem}
If $B$ is symmetric, we have that $\|A\otimes A\|_2=\|A\|_2^2=\lambda_{max}(A^TA)$ and that $\lambda_{max}(B+B^\top)=2\lambda_{max}(B)$. So we see that in this case our Condition \eqref{eq:assumptionpge1} implies Condition (4.4) of \cite[Theorem 4.5]{S2010}  for the existence of a stationary distribution for $p=k=1$.

 For a symmetric $B$ and $p=k>1$ the conditions are very similar only that an additional factor of $2^{p-1}$ appears inside the integral in the conditions ensuring geometric ergodicity. So the previously known conditions for the existence of a stationary distribution with a finite $p$-th moment for $p>1$ may be somewhat less restrictive than our conditions for geometric ergodicity with a finite $p$-th moment for $p>1$.  But in contrast to \cite{S2010} we do not  need to restrict ourselves to $B$ being diagonalizable and integer moments $p$. 
\end{bem}

	\begin{example}\label{condexamp1}
	Let us now consider that the driving L\'evy process is a compound Poisson process with rate $\gamma>0$ and jump distribution $P_L$. So $P_L$ is a probability measure on $\R^d$ and $\nu_L=\gamma P_L$.
	
	As in many applications where discrete time multivariate GARCH processes are employed the noise is taken to be an iid standard normal distribution, a particular choice would be to take $P_L$ as the standard normal law.
	
	However, we shall only assume that $P_L$ has finite second moments, mean zero and covariance matrix $\Sigma_L$. Observe that the upcoming Section \ref{section:irredresults} shows that the needed irreducibility and aperiodicity properties hold as soon as $P_L$ has an absolutely continuous component with a strictly positive density around zero (again definitely satisfied when choosing the standard normal distribution).
	
	Then we have that
	\begin{align}
	\left\|\int_{\R^{d}}  yy^\top  \nu_L(dy)\right\|_2&=\gamma\|\Sigma_L\|_2=\gamma\lambda_{max}(\Sigma_L),\\
	\int_{\R^{d}}\|y\|^2_2\nu_L(dy)&=\gamma\tr(\Sigma_L).
	\end{align}
	
	If we assume that $B=\beta I_d$, $A=\alpha I_d$ for some $\alpha,\beta \in \R$ and $\eta=I_d$, then \eqref{eq:cond1} and \eqref{eq:cond2} are equivalent to
	\[
	2\beta+\alpha^2\gamma \lambda_{max}(\Sigma_L)<0
	\]
	whereas 
	\eqref{eq:assumptionp<1}, \eqref{eq:assumptionpge1}, \eqref{eq:assumptionv} are becoming
	\[
		2\beta+\alpha^2\gamma \tr(\Sigma_L)<0
	\]
	for $p=1$. Note that in this particular set-up it is straightforward to see that the choice of $\eta$ has no effect on  \eqref{eq:cond1}, \eqref{eq:assumptionp<1}, \eqref{eq:assumptionpge1}, \eqref{eq:assumptionv}. The latter is also the case for \eqref{eq:cond2} if $\Sigma_L$ is additionally assumed to be a multiple of the identity. 
	
	As, for example, $\lambda_{max}(I_d)=1$, but $\tr(I_d)=d$, this also illustrates that Conditions 	\eqref{eq:assumptionp<1}, \eqref{eq:assumptionpge1}, \eqref{eq:assumptionv} are considerably more restrictive than Conditions \eqref{eq:cond1} and \eqref{eq:cond2} unless we are in the univariate case. 
	
	For $p\neq 1$  it is not sufficient to only specify the mean and variance of the jump distribution to check Conditions 	\eqref{eq:assumptionp<1}, \eqref{eq:assumptionpge1}, \eqref{eq:assumptionv}. For a concrete specification of $P_L$ it is, however, straightforward to check them by (numerical) integration.
		\end{example}
			\begin{example}\label{condexamp2}
				We consider the same basic set-up as in Example \ref{condexamp1} in dimension $d=2$.
				
				If we take $A=I_2,\,B=\begin{pmatrix}-2 & 0\\0 & -4\end{pmatrix},\,\gamma=1,\,\Sigma_L=\begin{pmatrix}3 & 0\\0 & 6\end{pmatrix} $ and $ \eta=I_2$, then Condition \eqref{eq:cond1} is violated whereas \eqref{eq:cond2} is satisfied.
				
				If we take $A=\begin{pmatrix}\sqrt3 & 0\\0 & \sqrt6\end{pmatrix} ,\,B=\begin{pmatrix}-2 & 0\\0 & -4\end{pmatrix},\,\gamma=1,\,\Sigma_L=I_2$ and $ \eta=I_2$, then Condition \eqref{eq:cond2} is violated whereas \eqref{eq:cond1} is satisfied.
				
				In both cases Conditions 	\eqref{eq:assumptionp<1}, \eqref{eq:assumptionpge1}, \eqref{eq:assumptionv} are violated for $p=1$.
			\end{example}
		\begin{remark}[General Matrix Subordinator]\label{rem:matsub}
			An inspection of our proofs shows, that we can also consider the stochastic differential equation
			\begin{equation}\label{eq:matsub}
			dY_{t} = (BY_{t-} + Y_{t-}B^{\top}) dt + A  (C+Y_{t-})^{\frac 1 2} dL_{t}(C+Y_{t-})_{t-}^{\frac 1 2} A^{\top}
			\end{equation}
			with  $A,B \in M_{d}(\R)$ and $C \in \s_{d}^{++}$, initial value $Y_{0} \in \s_{d}^{+}$ and $L$ being a $d\times d$ matrix subordinator with drift $\gamma_L=0$ and L\'evy measure $\nu_L$.
			
			All results we obtained in this section remain valid when replacing $\int_{\R^d}$ by $\int_{\sdd}$, $yy^T$ by $y$, $\|y\|_2^2$ by $\|y\|_2$ and $\E(\|L_1\|_2^{2p}) < \infty$ by $\E(\|L_1\|_2^{p}) < \infty$.
			
		\end{remark}
		\begin{remark}[MUCOGARCH(p,q)]
			In principle, it is to be expected that the results can be generalized with substantial efforts to general order MUCOGARCH(p,q) processes (cf. \cite{fuchs,DettlingStelzer2019} for discrete time BEKK models). However, as explained in the introduction, they have only been defined briefly and no applicable conditions on the admissible parameters are known nor is the natural state space of the associated Markov process (which should not be some power of the positive definite cone).  As can be seen from the discrete time analogue, the main issue is to find a suitable Lyapunov test function and suitable weights need to be introduced. For orders different from (1,1) the weights known to work in the discrete time setting are certainly not really satisfying except when the Lyapunov test function is related to the first moment (the case analysed in \cite{fuchs}).  On the other hand in applications very often order (1,1) processes provide rich enough dynamics already. Therefore, we refrain from pursuing this any further.
		\end{remark}
		
		The geometric ergodicity results of this section are of particular relevance for at least  the following:
		\begin{enumerate}
			\item Model choice:\\
			In many applications  stationary models are called for and one should use models which have a unique stationary distribution. Our results are the first giving sufficient criteria for a MUCOGARCH model to have a unique stationary distribution (for the volatility process).
			\item Statistical inference:\\
			\cite{doRegoSousaKlueppelbergStelzer2018, doRegoSousa2019} investigate  in detail a moment based estimation method for the parameters of MUCOGARCH(1,1) processes. This involves establishing identifiability criteria, which is a highly non-trivial issue, and calculating moments explicitly, which requires assumptions on the moments of the driving L\'evy process. However, the asymptotics of the estimators derived there hinge  centrally on the results of the present paper. Actually, we cannot see any reasonable other way to establish consistency and asymptotic normality of estimators for the MUCOGARCH parameters  than using the geometric ergodicity conditions for the volatility of the present paper and to establish that under stationarity it implies that the increments of the MUCOGARCH process $G$ are ergodic and strongly mixing (see 	\cite{doRegoSousaKlueppelbergStelzer2018, doRegoSousa2019} for details).
			\item Simulations:\\
			Geometric ergodicity implies that simulations (of a Markov process)  can be started with an arbitrary initial value and that after a (not too long) burn-in period the simulated path behaves essentially like one following the stationary dynamics. On top $V$-uniform ergodicity (which is actually obtained above) provides results on the finiteness of certain moments and the convergence of the moments to the stationary case. We illustrate this now in some exemplary simulations.
		\end{enumerate}
		\begin{example}\label{ex:sim}
					We consider the set-up of Example \ref{condexamp1} in dimension two and choose $\alpha=0.14, \beta=-0.01, C=I_2, \gamma=1$ and the jumps to be standard normally distributed.  
					
			\begin{figure}
				\begin{center}
					\includegraphics[width=0.49\linewidth]{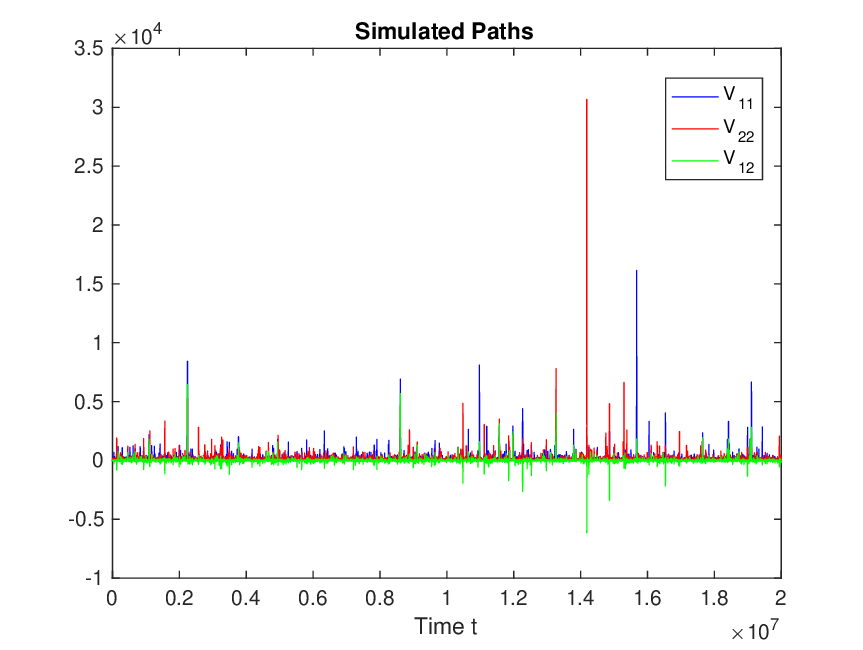} \includegraphics[width=0.49\linewidth]{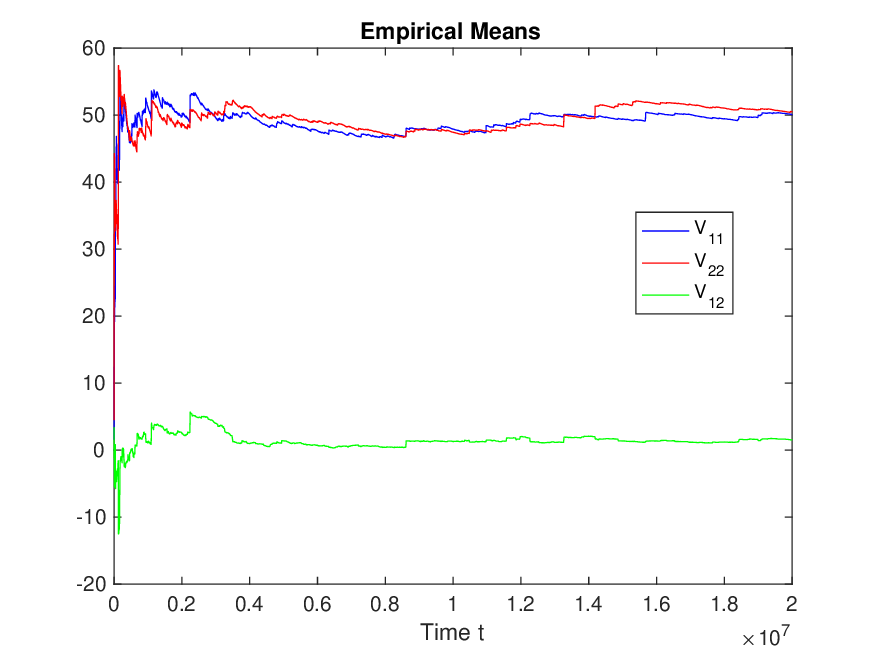}\\
					\includegraphics[width=0.49\linewidth]{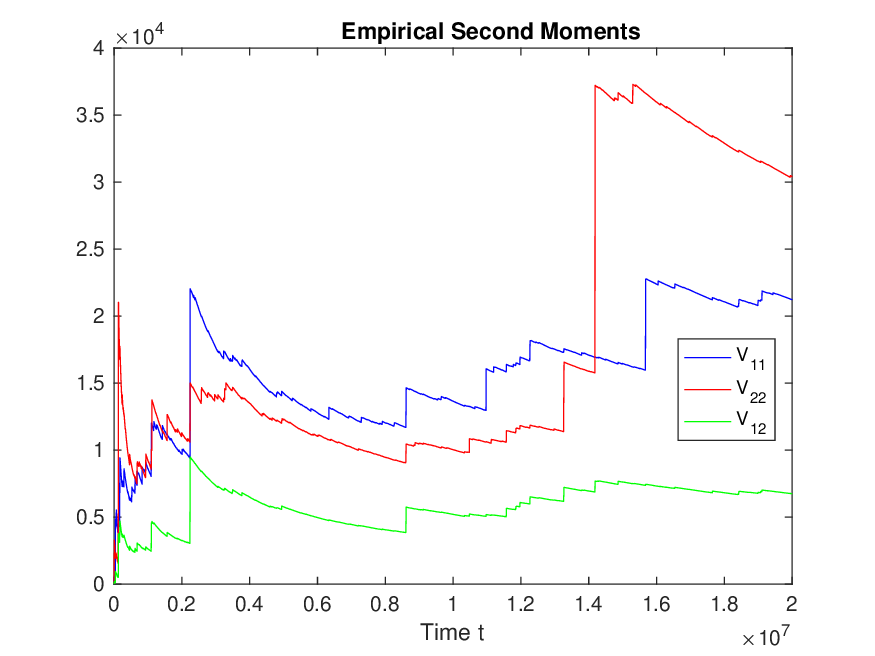} \includegraphics[width=0.49\linewidth]{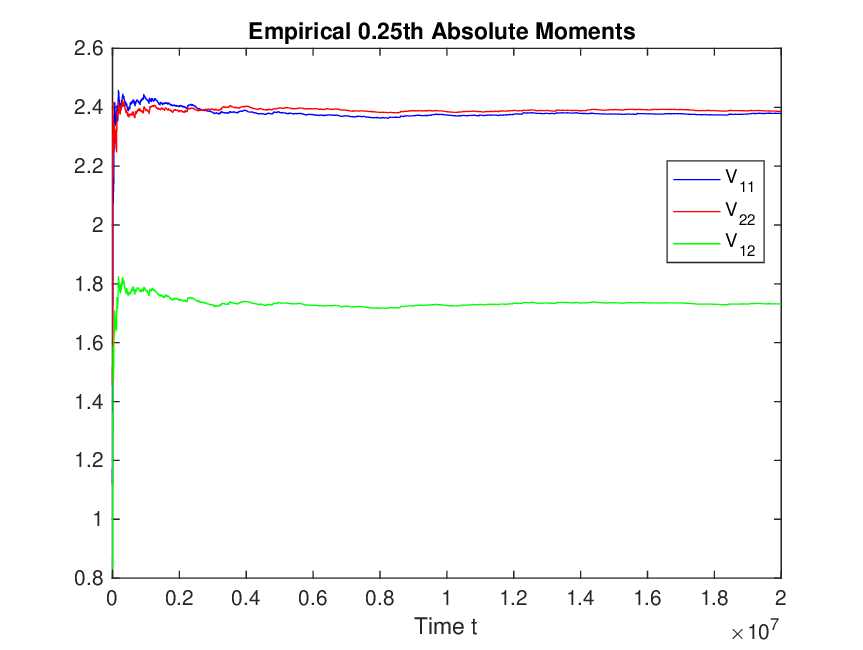}
				\end{center}
				\caption{Plots of the simulations of Example \ref{ex:sim} for $\alpha=0.14$. The first variance component ($V_{11}$) is depicted in blue, the second ($V_{22}$) in red and the covariance component ($V_{12}$) in green. The time  horizon is zero to 20 million.}\label{fig:140}
			\end{figure}
			\begin{figure}
				\begin{center}
					\includegraphics[width=0.49\linewidth]{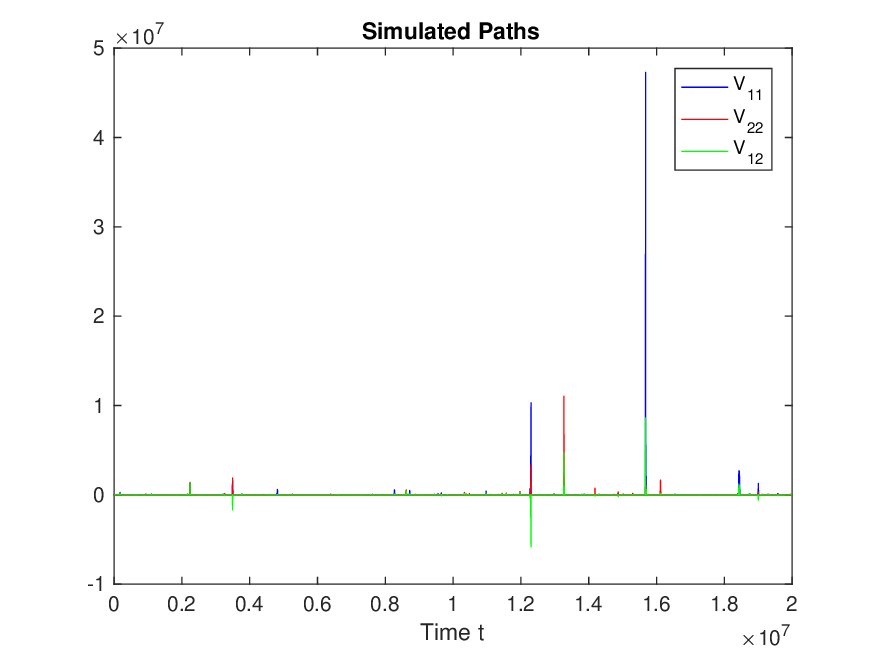} \includegraphics[width=0.49\linewidth]{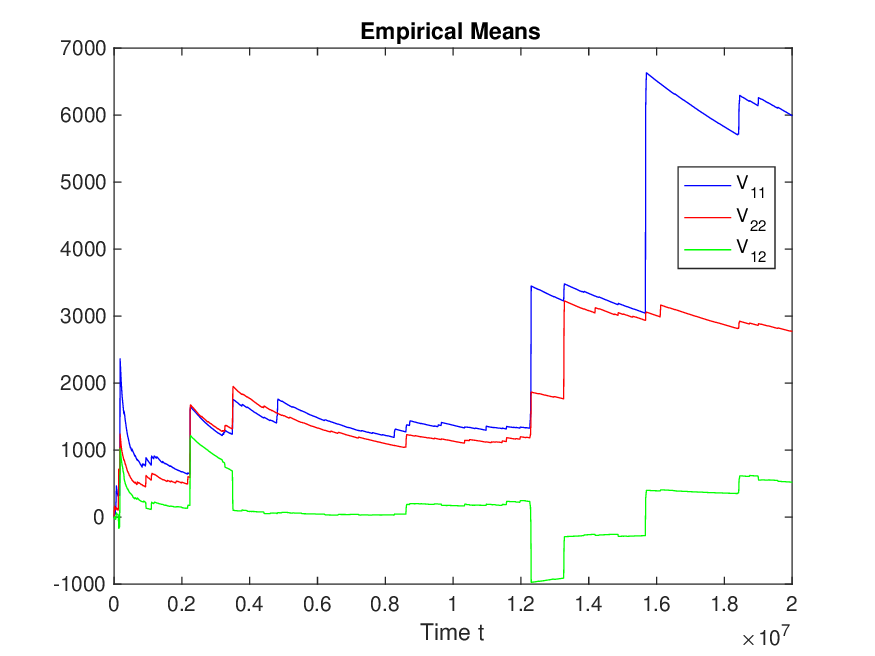}\\
					\includegraphics[width=0.49\linewidth]{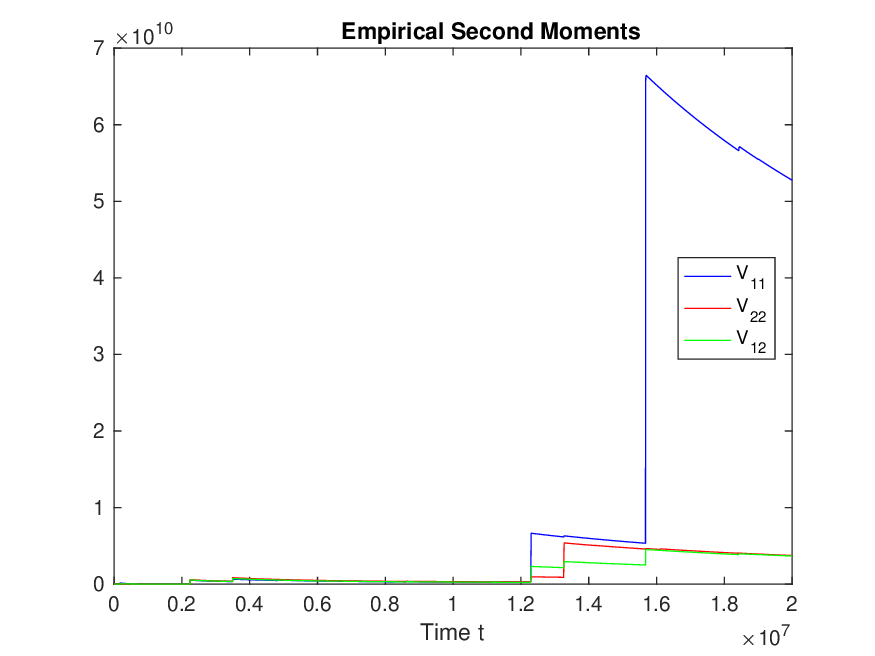} \includegraphics[width=0.49\linewidth]{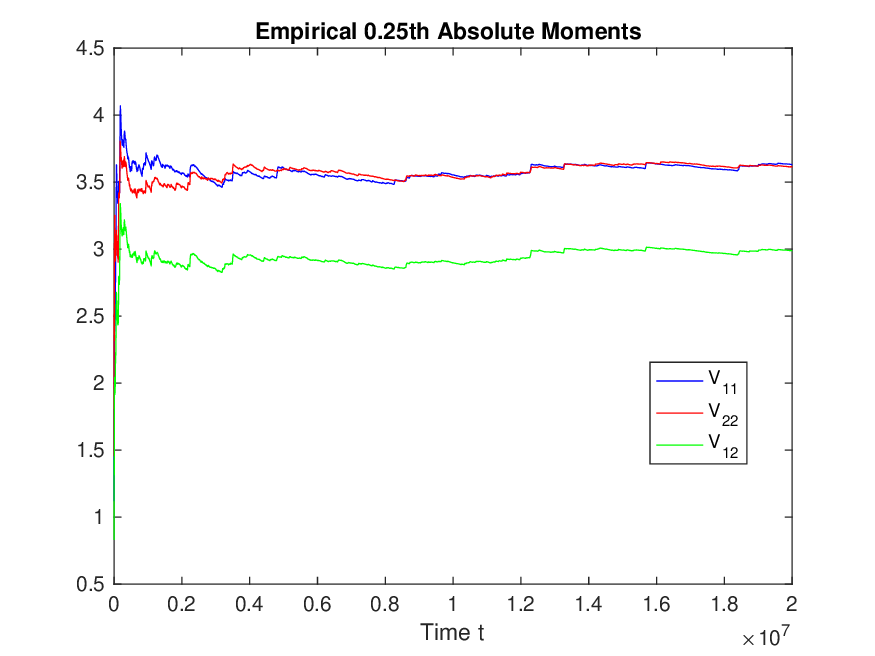}
				\end{center}
				\caption{Plots of the simulations of Example \ref{ex:sim} for $\alpha=0.142$. The first variance component ($V_{11}$) is depicted in blue, the second ($V_{22}$) in red and the covariance component ($V_{12}$) in green. The time  horizon is zero to 20 million.}\label{fig:142}
			\end{figure}
			\begin{figure}
				\begin{center}
					\includegraphics[width=0.49\linewidth]{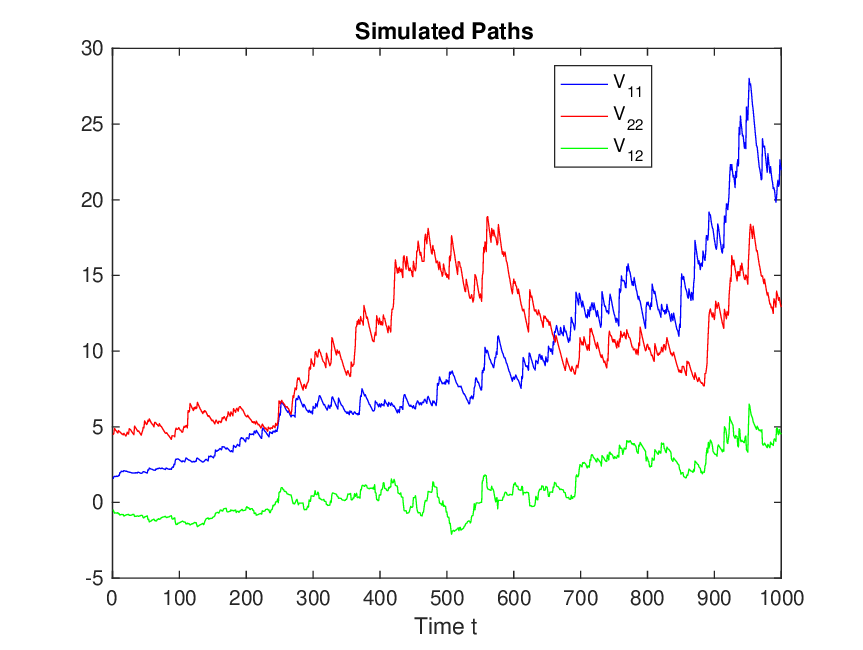} 		\includegraphics[width=0.49\linewidth]{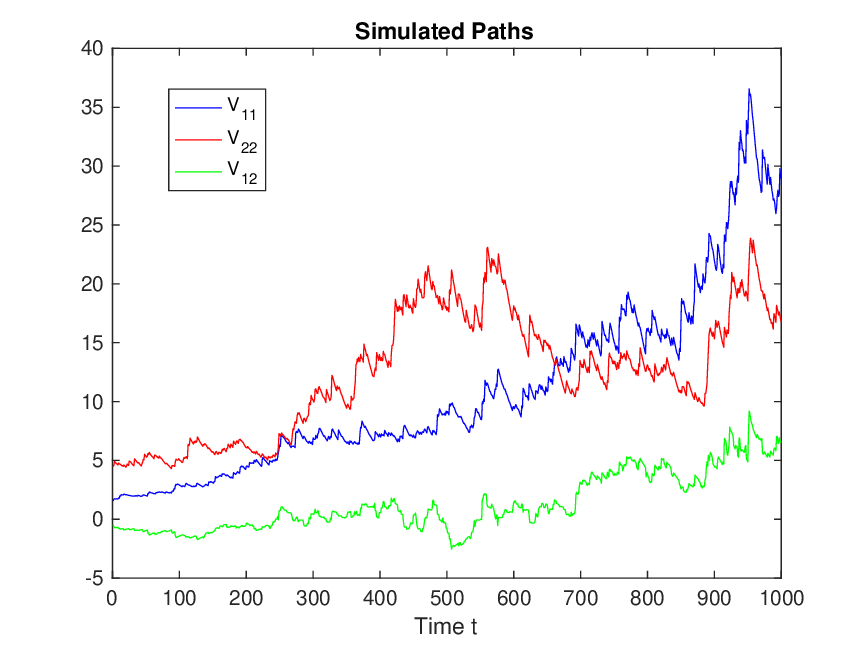} 
				\end{center}
				\caption{Plots of the simulations of Example \ref{ex:sim}, $\alpha=0.14$ on the left and $\alpha=0.142$ on the right. The first variance component ($V_{11}$) is depicted in blue, the second ($V_{22}$) in red and the covariance component ($V_{12}$) in green. The time  horizon is zero to 1000.}\label{fig:initpaths}
			\end{figure}
	
			Then condition \eqref{eq:cond1} (or \eqref{eq:cond2}) is satisfied (but not conditions \eqref{eq:assumptionp<1}, \eqref{eq:assumptionpge1}, \eqref{eq:assumptionv}) and so we have geometric ergodicity with convergence of the (absolute) first moment. However, this is a rather extreme case as (given all the other parameter choices)  \eqref{eq:cond1} is satisfied if and only if $|\alpha|<\sqrt {0.02}\approx 0.1414$ and thus some convergences are only  to be seen clearly for very long paths. We have simulated a path up to time 20 million with Matlab (see Figure \ref{fig:140}). The paths themselves look immediately stationary (at this time scale) but extremely heavy tailed. In the long run (and not too fast) the running empirical means are reasonably converging (to their true values of 50 for the variances and 0 for the covariance). In sharp contrast to this there seems to be no convergence for the empirical second moments at all and actually we strongly conjecture that the stationary covariance matrix is infinite in this case. Looking at lower moments, like the 0.25th (absolute) empirical moments there is again convergence to be seen which is faster than the one for the first moment (as is to be expected). Zooming into the paths at the beginning in Figure \ref{fig:initpaths} (left plot)  shows that the process clearly does not start of like a stationary one, but the behaviour looks stationary to the eye very soon.	
			
			Changing only the parameter $\alpha$ to $0.142$ and using the same random numbers we obtain the results depicted in Figure \ref{fig:142}. The initial part of the paths (see  right plot in Figure \ref{fig:initpaths}) looks almost unchanged (observe that the scaling of the vertical axes in the plots of Figure  \ref{fig:142} differs). The paths in total clearly appear to be even more heavy tailed. Now neither the empirical second moments nor the empirical means seem to converge (and we conjecture  that both mean and variance are infinite for the stationary distribution) . Now condition \eqref{eq:cond1}  does not hold any more and so the simulations suggest that this condition implying geometric ergodicity  with convergence of the first moments is pretty sharp, as expected from the theoretical insight. The simulations also seem to clearly indicate that the 1/4th (absolute) moments still converge. Numerical integration shows for  \eqref{eq:assumptionp<1} (again the choice of $\eta$ makes no difference in this special case):
				\begin{equation*}\int_{\R^{2}} \left( \left( 1+\alpha^2\|y\|^2_2\right)^\frac14 - 1\right)\nu_L(dy)
				+2\beta/4\approx 0.0098  -0.005   > 0.
				\end{equation*}
				So our sufficient condition  \eqref{eq:assumptionp<1} for convergence of the $1/4$th moment is not satisfied, but due to the involved inequalities we expected it not to be too sharp.  Actually, numerically not even the logarithmic condition \cite[(4.3)]{S2010} for the existence of a stationary distribution seems to be satisfied. This illustrates that these conditions involving norm estimates are not too sharp, as our simulations and the fact that condition   \eqref{eq:cond1} is almost satisfied strongly indicate that geometric ergodicity with convergence of some moments $0<p<1$ should hold.
		\end{example}
		Finally, we consider an example illustrating that our conditions for $p\not =1$ are more precise if one coordinate dominates the others.
		\begin{example}\label{ex:sim2}
			We consider the same model as in the previous Example \ref{ex:sim} with the only difference that we now assume the jumps to have a two-dimensional normal distribution with covariance matrix $\Sigma_L=\begin{pmatrix}
			1 & 0\\ 0 & \sigma
			\end{pmatrix}$ for $\sigma>0$. For $\sigma\to 0$ we conclude from the formulae given in Example \ref{condexamp1} that for $p=1$ \eqref{eq:assumptionp<1}, \eqref{eq:assumptionpge1}, \eqref{eq:assumptionv} are asymptotically equivalent to \eqref{eq:cond1} and as in Example  \ref{ex:sim} the choice of $\eta$ does not matter (they also are equivalent to \eqref{eq:cond2}  for $\eta=I_2$, but there the choice of $\eta$ may now matter). Note that  \eqref{eq:cond1} is not affected by changes in $\sigma$ as long as $\sigma\leq 1$.
			 
			 Choosing $\sigma=1/1000$ and $\alpha=0.142$ as in the second part of Example  \ref{ex:sim} numerical integration shows for  \eqref{eq:assumptionp<1}  with $p=1/4$:
			 \begin{equation*}\int_{\R^{2}} \left( \left( 1+\alpha^2\|y\|^2_2\right)^\frac14 - 1\right)\nu_L(dy)
			 +2\beta/4\approx 0.0049  -0.005   < 0.
			 \end{equation*}
			 So we obtain geometric ergodicity and the finiteness of the $1/4$th absolute moments from our theory.
			 
			 	\begin{figure}
			 		\begin{center}
			 			\includegraphics[width=0.49\linewidth]{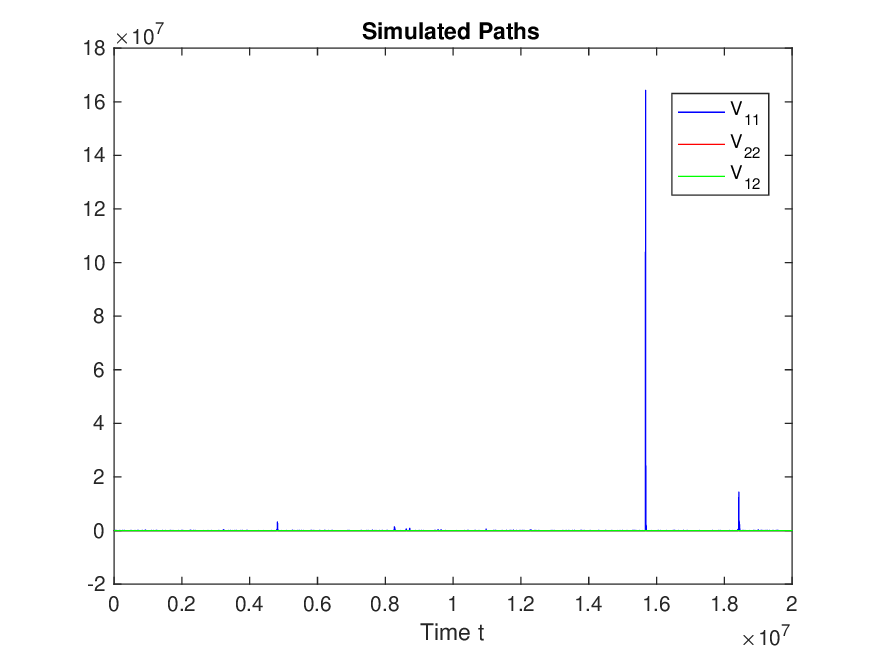} \includegraphics[width=0.49\linewidth]{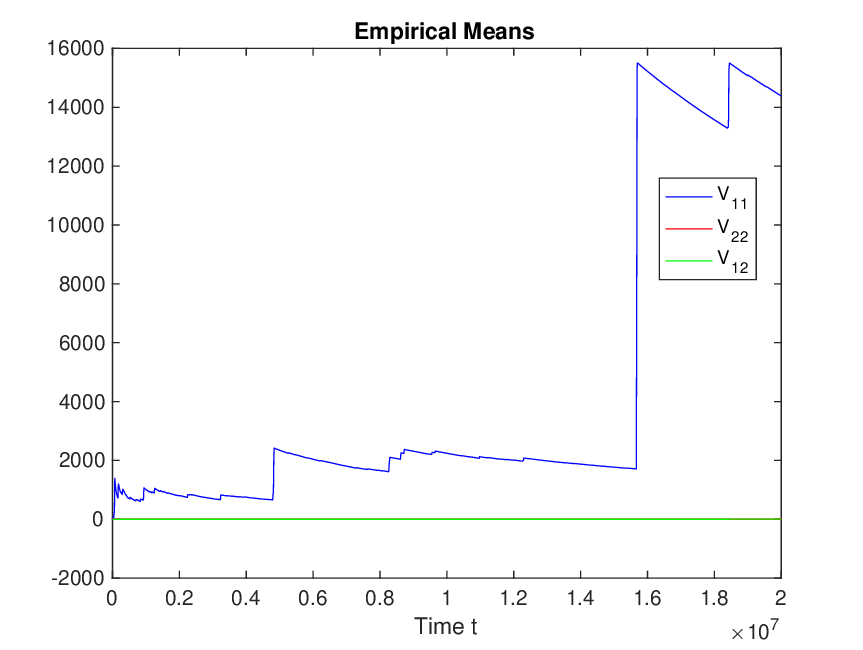}\\
			 			\includegraphics[width=0.49\linewidth]{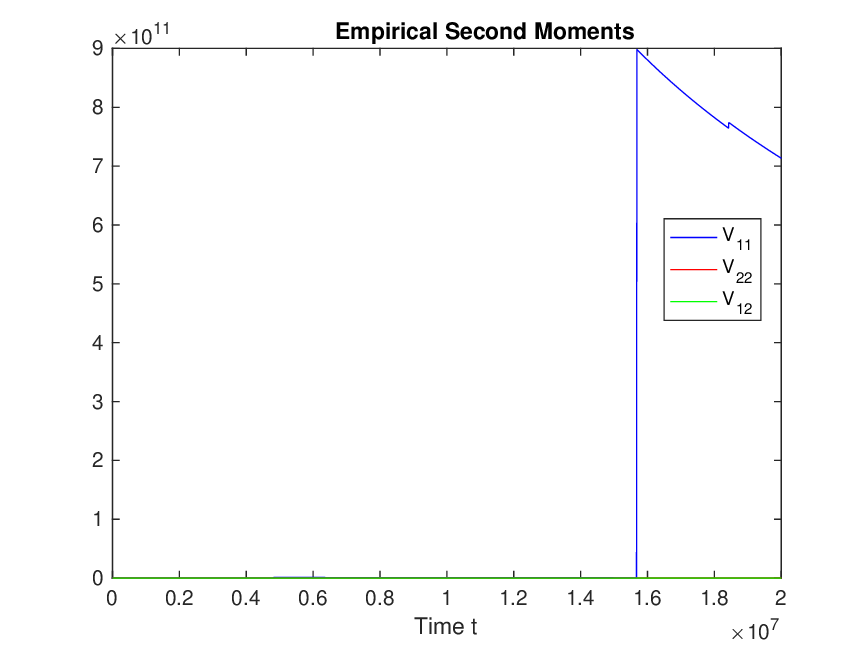} \includegraphics[width=0.49\linewidth]{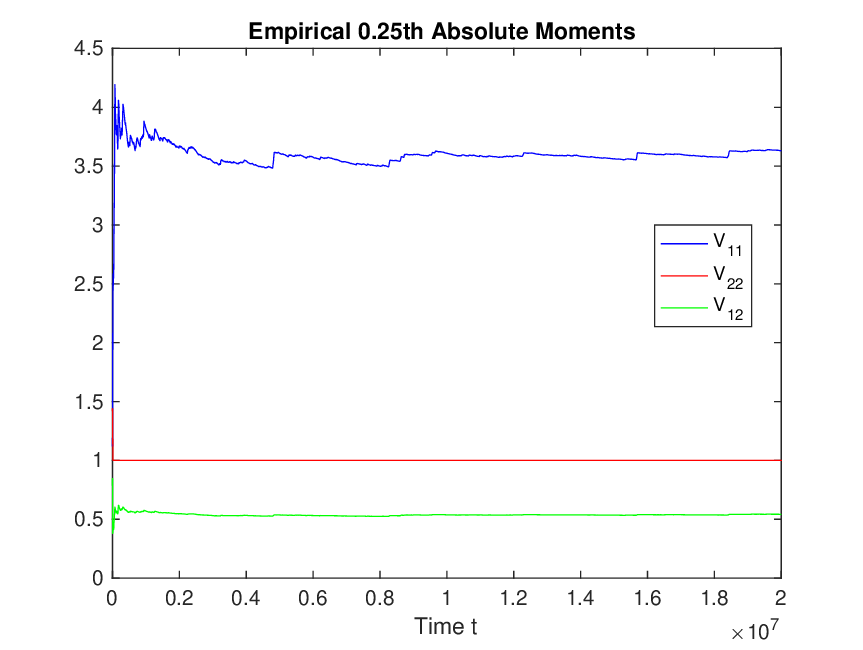}
			 		\end{center}
			 		\caption{Plots of the simulations of Example \ref{ex:sim2}. The first variance component ($V_{11}$) is depicted in blue, the second ($V_{22}$) in red and the covariance component ($V_{12}$) in green. The time  horizon is zero to 20 million.}\label{fig:nd}
			 	\end{figure}
			 	
			 	Figure \ref{fig:nd} shows the results of a simulation using the same time horizon and random numbers as in Example \ref{ex:sim}. The plots of the paths, the empirical means and second moments) show that the first variance component extremely dominates the others. One immediately sees that the empirical mean and second moment of $V_{11}$ does not converge. In line with our theory all 0.25th absolute moments seem to converge pretty fast.
		\end{example}
\section{Irreducibility and aperiodicity of the MUCOGARCH volatility process} 
\label{section:irredresults}
So far we have assumed away the issue of irreducibility and aperiodicity focusing on establishing a Foster-Lyapunov drift condition to obtain geometric ergodicity.  As this is intrinsically related to the question of the existence of transition densities and support theorems, this is a very hard problem of its own. In our case the degeneracy of the noise, i.e. that all jumps of the L\'evy process are rank one and that thus the jump distribution has no absolutely continuous component, makes our case particularly challenging.

We now establish sufficient conditions for the irreducibility and aperiodicity of $Y$ by combining at least $d$ jumps to get an absolutely continuous component of the increments. 

\begin{theorem}[Irreducibility and Aperiodicity] \label{theo:irreducible} 
Let $Y$ be a MUCOGARCH volatility process driven by a compound 
Poisson process $L$ and with $A\in GL_d(\R)$, $\Re(\sigma(B))<0$. If the jump distribution of $L$ has a non-trivial absolutely 
continuous component equivalent to the Lebesgue measure on $\R^d$, 
then  $Y$ is irreducible with respect to the Lebesgue measure on  
$\s_d^+$ and aperiodic. 
\end{theorem}

We can soften the conditions on the jump distribution of the Compound Poisson process:

\begin{cor} \label{cor:irredcppnullenvironment} Let $Y$ be a MUCOGARCH volatility process driven by a compound Poisson process $L$ and with $A\in GL_d(\R)$, $\Re(\sigma(B))<0$.  If the jump distribution of $L$ has a non-trivial absolutely continuous component  equivalent to the Lebesgue measure on $\R^d$ restricted to an open neighborhood of zero, then $Y$ is irreducible w.r.t. the Lebesgue measure restricted to an open neighborhood of zero in $\s_d^+$ and aperiodic.
\end{cor}

\begin{bem}
(i) If the driving L\'evy process is an infinite activity L\'evy process, whose L\'evy measure has a non-trivial absolutely continuous component and the density of the component w.r.t. the Lebesgue measure on $\s_d^+$ is strictly positive in a neighborhood of zero, we can show that $Y$ is open-set irreducible w.r.t. the Lebesgue measure restricted to an open neighborhood of zero in $\sdd$. For strong Feller processes open-set irreducibility provides irreducibility, but the strong Feller property is to the best of our knowledge hard to establish for L\'evy-driven SDEs. A classical way to show irreducibility is by using density or support theorems based on  Malliavin Calculus, see e.g. \cite{denis2000, kunita2009, simon2000}. But they all require, that the coefficients of the SDE have bounded derivatives, which is not the case for the MUCOGARCH volatility process. So finding criteria for irreducibility in the infinite activity case appears to be a very challenging question beyond the scope of the present paper. 

(ii) In the univariate case establishing irreducibility and aperiodicity is much easier, as one can easily use the explicit representation of the COGARCH(1,1) volatility process to establish that one has weak convergence to a unique stationary distribution which is self-decomposable and thus has a strictly positive density. This way is blocked in the multivariate case as we do not have the explicit representation and the linear recurrence equation structure.
\end{bem}

Note that the condition that the L\'evy measure has an absolutely continuous component with a support containing zero is the obvious analogue to the condition on the noise in \cite{fuchs}.

\begin{remark}
Let us return to the SDE \eqref{eq:matsub} driven by a general matrix subordinator which we introduced in Remark \ref{rem:matsub}. It is straightforward to see that then the conclusions of Theorem \ref{theo:irreducible} or Corollary \ref{cor:irredcppnullenvironment}, respectively, are valid for the process $Y$, if the matrix subordinator $L$ is a compound Poisson process with the jump distribution of $L$ having a non-trivial absolutely 
continuous component equivalent to the Lebesgue measure on $\s_d^+$   (restricted to an open neighborhood of zero).
\end{remark}

\section{Proofs}\label{section:proofsMarkov}
To prove our results we use the stability concepts for Markov processes of \cite{MT2, MT3}.
\subsection{Markov processes and ergodicity} \label{chaptermarkov}
For the convenience of the reader, we first give a short introduction to the definitions and results for general continuous time Markov processes following mainly
 \cite[Section 3]{DMT1995}.

Given be  a state space $X$, which is assumed to be an open or closed subset of a finite-dimensional real vector spaced equipped with the usual topology and the Borel-$\sigma$-algebra $\B(X)$. We consider a continuous time Markov process $\Phi=(\Phi_t)_{t\ge0}$  on $X$   with transition probabilities 
$\p^t(x,A)=\p_x(\Phi_t \in A)$ for $x\in X, A \in \B(X)$. 

The operator $\p^t$ from the associated transition semigroup acts on a bounded 
measurable function $f$ as
\begin{equation}
\p^t f(x) = \int_X \p^t(x,dy)f(y)
\end{equation}
and on a $\sigma$-finite measure $\mu$ on $X$ as
\begin{equation}
\mu \p^t(A) = \int_X \mu(dy) \p^t(y,A).
\end{equation}

To define non-explosivity, we consider a fixed family $\{ O_n | n \in \Z_+\}$ 
of open precompact sets, i.e. the closure of $O_n$ is a compact subset of $X$, 
for which $O_n \nearrow X$ as $n \to \infty$. With $T^m$ we denote the first 
entrance time to $O_m^c$ and by $\xi$ the exit time for the process, defined as
$
 \xi \overset{\mathcal D}{=} \lim_{m \to \infty} T^m.
$

\begin{definition}[Non-explosivity, {\cite[Chapter 1.2]{MT3}}]
 We call the process $\Phi$ \textbf{non-explosive} if $\p_x(\xi = \infty)=1$ 
for all $x\in X$.
\end{definition}

 By $\cb(X)$ 
we denote the set of all continuous and bounded functions $f: X \to \R$ and 
by $\co(X)$ those continuous and bounded functions, which vanish at infinity. 

\begin{definition}[Stochastic continuity and Feller processes]
	Let $(\p^t)_{t\in\R_+}$ be the transition semigroup of a time homogeneous 
	Markov process $\Phi$. \\
	(i) $(\p^t)_{t\in\R_+}$ or $\Phi$ is called 
	\textbf{stochastically continuous} if 
	\begin{equation}
	\lim_{\substack{t\to 0 \\ t\ge 0}} \p^t(x, \mathcal N(x)) = 1
	\end{equation}
	for all $x \in X$ and open neighborhoods $\mathcal N(x)$ of $x$. \\
	(ii) $(\p^t)_{t\in\R_+}$ or $\Phi$ is a \textbf{(weak) 
		$\mathbf{\cb}$-Feller} semigroup or process if it is stochastically continuous 
	and 
	\begin{equation} \label{cbfeller}
	\p^t (\cb(X)) \subseteq \cb(X) \text{  for all  } t\ge 0. \\
	\end{equation}
	(iii) If in (ii) we have instead of (\ref{cbfeller})
	\begin{equation}
	\p^t (\co(X)) \subseteq \co(X) \text{  for all  } t\ge 0
	\end{equation}
	we call the semigroup or the process \textbf{(weak) $\mathbf{\co}$-Feller}.
\end{definition}
 Combining the definition of strongly continuous contraction 
 semigroups, Theorem 4.1.1  and the definition of Feller processes in  \cite{marcusrosen}  shows that a $ {\co}$-Feller process is a Borel-right process (cf.  \cite{marcusrosen} for a definition).

From now on we  assume that $\Phi$ is a non-explosive Borel right process. For the definitions and details of 
the existence and structure see \cite{sharpe}.

\begin{definition}[Invariant measure, {\cite[Chapter 3 ]{DMT1995}}]
A $\sigma$-finite measure $\pi$ on $\B(X)$ with the property 
\begin{equation}
\pi = \pi \p^t, \  \forall t \ge 0                                                              
\end{equation}
is called \textbf{invariant}.
 
\end{definition}

Notation: By $\pi$ we always denote an invariant measure of $\Phi$, if it exists. 

\begin{definition}[Exponential ergodicity, {\cite[Chapter 3 ]{DMT1995}}]
 $\Phi$ is called \textbf{exponentially ergodic}, 
if an invariant measure $\pi$ exists and satisfies for all $x \in X$
\begin {equation}
 \label{expergodic} \| \p^t(x,.) - \pi\|_{\text{TV}} \le M(x) \rho^t, \  \  \forall t\ge0
\end {equation}
for some finite $M(x)$, some $\rho<1$ and where $\|\mu\|_{\text{TV}}:=\sup_{|g|\le 1,\,g \text{ measurable}} | \int \mu(dy)g(y)|$  denotes the total variation norm. 

If this convergence holds for the $f$-norm $\|\mu\|_f:=\sup_{|g|\le f,\,g \text{ measurable}} | \int \mu(dy)g(y)|$ (for any signed measure $\mu$) , 
where $f$ is a measurable function from the state space $X$ to $[1,\infty)$, we 
call the process \textbf{$\mathbf f$-exponentially ergodic}.
\end{definition}

A seemingly stronger formulation of $V$-exponential ergodicity is $V$-uniform ergodicity: We require that $M(x) = V(x) \cdot D$ with some finite constant $D$. 
\begin{definition}[$V$-uniform ergodicity, {\cite[Chapter 3]{DMT1995}}]
 $\Phi$ is called \textbf{$\mathbf V$-uniformly ergodic}, if a measurable 
function $V: X \to [1,\infty)$ exists such that for all $x \in X$ 
\begin{equation}
 \| \p^t(x,.) - \pi \|_V \le V(x) D \rho^t, \  t\ge0
\end{equation}
holds for some $D<\infty, \rho <1$.
\end{definition}

To prove ergodicity we need the notions of irreducibility and aperiodicity. 

\begin{definition}[{\cite{DMT1995}}, Chapter 3]
 For any $\sigma$-finite measure $\mu$ on $\mathcal B(X)$ we call the process 
$\Phi$ {\boldmath$ \mu$}\textbf{-irreducible} if for any $B \in  \mathcal B(X)$ 
with $\mu(B)>0$ 
\begin{equation}
 \E_x (\eta_B) > 0, \forall x \in X
\end{equation}
holds, where $\eta_B := \int_0^\infty \1_{\{\Phi_t \in B\} } dt$ is the occupation time.
\end{definition}

This is obviously the same as requiring 
\begin{equation}
 \int_0^\infty \p^t(x,B) dt >0,\  \forall x \in X.
\end{equation}

If $\Phi$ is $\mu$-irreducible, there exists a maximal irreducibility measure $\psi$ 
such that every other irreducibility measure $\nu$ is absolutely continuous with respect to $\psi$ (see \cite[p. 493]{MT2}). 
We write $\boldsymbol{\mathcal B^+(X)}$ for the collection of all measurable subsets $A\in \mathcal B(X)$ with $\psi(A)>0$. 

\begin{bem} \label{irredskeleton}
In \cite[Proposition 1.1]{TT1979} it was shown, that if the discrete time $h$-skeleton of a process, the $\p^h$-chain, 
is $\psi$-irreducible for some $h>0$, 
then it holds for the continuous time process. If the $\p^h$-chain  is $\psi$-irreducible for every $h>0$, 
we call the process \textbf{simultaneously $\boldsymbol{\psi}$-irreducible}.
\end{bem}

One probabilistic form of stability is the concept of Harris recurrence.

\begin{definition}[Harris recurrence, {\cite[Chapter 2.2]{MT2}}]\ \\
(i)\ $\Phi$ is called \textbf{Harris recurrent}, if either
		\begin{itemize}
		\item  $\p_x(\eta_A = \infty)=1$ whenever $\phi(A)>0$ for some $\sigma$-finite measure $\phi$, or
		\item $\p_x(\tau_A < \infty)=1$ whenever $\mu(A)>0$ for some $\sigma$-finite measure $\mu$. $\tau_A:= \inf \{t\ge0~:~ \Phi_t \in A \}$ is the first hitting time of $A$. 
	\end{itemize}
	(ii)\ Suppose that $\Phi$ is Harris recurrent with finite invariant measure $\pi$, then $\Phi$ is called \textbf{positive Harris recurrent}.
\end{definition}

To define the class of subsets of $X$ called petite sets, we suppose that $a$ is a probability distribution 
on $\R_+$. We define the Markov transition function $K_a$ for the process sampled by $a$ as
\begin{equation}
 K_a(x,A):= \int_0^\infty \p^t(x,A)a(dt), \  \forall x\in X, A\in\mathcal B.
\end{equation}
\begin{definition}[Petite and small sets, {\cite[Chapter 3]{DMT1995}}]
 A nonempty set $C\in \mathcal B$ is called \textbf{$\boldsymbol{\nu_a}$-petite}, if $\nu_a$ is a nontrivial measure on $\mathcal B(X)$ 
and $a$ is a sampling distribution on $(0,\infty)$ satisfying
\begin{equation}
 K_a(x,.) \ge \nu_a(.), \  \forall x\in C.
\end{equation}
  When the sampling distribution $a$ is degenerate, i.e. a single point mass, we call the set $C$ \textbf{small}. 
\end{definition}

\begin{bem} \label{remark:small} 
 Like in the discrete time Markov chain theory the set $C$ is small, if there exists an $m>0$ and a nontrivial measure
  $\nu_m$ on $\mathcal B(X)$ such that for all $x\in C,B\in \mathcal B(X)$
 \begin{equation}
  \p^m(x,B) \ge \nu_m(B)
 \end{equation}
 holds. 
\end{bem}

For discrete time chains there exists a well known concept of periodicity, see for example \cite[Chapter 5.4]{MT}.
For continuous time processes this definition is not adaptable, since there are no fixed time steps. 
But a similar concept is the definition of aperiodicity for continuous time Markov processes as introduced in \cite{DMT1995}.

\begin{definition}[{\cite{DMT1995}}, Chapter 3] \label{def:aperiodic}
 A $\psi$-irreducible Markov process is called \textbf{aperiodic} if for some small set $C \in \mathcal B^+(X)$ there exists a $T$ 
 such that 
 $\p^t(x,C) >0$ for all $t\ge T$ and all $x\in C$. 
\end{definition}

 \begin{prop} \label{prop:simirred}
When $\Phi$ is simultaneously $\psi$-irreducible then we know from \cite[Proposition 1.2]{TT1979} that
 every skeleton chain is aperiodic in the sense of a discrete time Markov chain.
 \end{prop}
 
 For discrete time Markov processes there exist conditions such that every compact set is petite and every petite set is small:
 \begin{prop}[{\cite{MT}}, Theorem 6.0.1 and 5.5.7] \label{prop:compactpetite} \ \\ 
(i)\ If $\Phi$, a discrete time Markov process, is a $\Psi$-irreducible Feller chain with $supp(\Psi)$ having non-empty interior, every compact set is petite. \\
(ii) \ If $\Phi$ is irreducible and aperiodic, then every petite set is small.  
 \end{prop}
\begin{bem}
 Proposition \ref{prop:compactpetite} (i) is also true for continuous time 
Markov processes, see \cite{MT2}.
\end{bem}

To introduce the Foster-Lyapunov criterion for ergodicity we need the concept of the extended generator of a Markov process.

\begin{definition}[Extended generator, {\cite[Chapter 4]{DMT1995}}]
$\mathcal D (\A)$ denotes the set of all functions $f: X \times \R_{+} \to \R$ 
for which a function $g: X \times \R_{+} \to \R$ exists, such that $\forall x 
\in X, t>0$
\begin{align}
\label{exgen1} &\E_{x}(f(\Phi_{t},t)) = f(x,0) + \E_{x}\left( \int_{0}^{t}g(\Phi_{s},s)ds\right), \\
 \label{exgen2} &\E_{x}\left( \int_{0}^{t}|g(\Phi_{s},s)|ds\right) < \infty
\end{align}
holds. We write {\boldmath$\A f := g$} and call $\A$ the \textbf{extended generator of }{\boldmath $\Phi$}. 
$\mathcal D(\A)$ is called the domain of $\A$. 
\end{definition}
 
The next theorem from \cite{DMT1995} gives for an irreducible and aperiodic 
Markov process a sufficient criterion to be
$V$-uniformly ergodic. This is a modification of the Foster-Lyapunov drift 
criterion of \cite{MT3}. 

\begin{theorem}[{\cite{DMT1995}}, Theorem 5.2] \label{dmtergodicity}
 Let $(\Phi_t)_{t\ge 0}$ be a $\mu$-irreducible and aperiodic Markov process. If 
there exist 
constants $b,c > 0$ and a petite set $C$ in $\B(X)$ as well as a measurable function $V: 
X \to [1,\infty)$ such that \begin{equation}
 \A V \le -bV + c \1_C ,
\end{equation}
where $\A$ is the extended generator, then $(\Phi_t)_{t\ge 0}$ is 
$V$-uniformly ergodic.
\end{theorem}

\subsection{Proofs of Section \ref{section:mainresults}} 
\label{subsection:proofs}
We prove the geometric ergodicity of the MUCOGARCH volatility process by using 
Theorem \ref{dmtergodicity}. The main task is to show the validity of a Foster-Lyapunov drift condition and that the used function belongs to the domain of the extended generator. For the latter we need the existence of moments given in Lemma \ref{lemma:moments}. The different steps require similar inequalities. The most precise estimations are needed for the Foster-Lyapunov drift condition. Therefore we below first prove the Forster-Lyapunov drift condition assuming for a moment that we are in the domain of the extended generator and prove this and the finiteness of moments later on, as the inequalities needed there are often obvious from the proof of the Foster-Lyapunov drift condition.

\subsubsection{Proof of Theorem \ref{theoremge}}

The first and main part of the proof is to show the geometric ergodicity, as the positive Harris recurrence is essentially  a consequence of it.  
To prove the geometric ergodicity, and hence the existence and uniqueness of a stationary distribution, of the MUCOGARCH volatility process, it is enough to show that the Foster-Lyapunov drift condition of Theorem \ref{dmtergodicity} holds. All other conditions of Theorem \ref{dmtergodicity} follow from Theorem \ref{th:markov} or  the assumptions.


\paragraph{The extended generator} ~\\
Using e.g. the results on SDEs, the stochastic symbol and the infinitesimal generator of  \cite[Chapter 6.1]{schnurr2009}, \cite[Chapter 6.1]{schnurr2009} and \cite[Theorem 3.1]{SS2010} and that $\tr(XY)$ is the canonical scalar product in $\s_d$, we see that for a sufficiently regular function $u: \s_d^+\to\R$ in the domain of the (extended/infinitesimal) generator $\mathcal{A}$ we have.
\begin{align*}
 \A u(x) &= \tr\left((Bx+xB^\top) \nabla u(x)\right) + \int_{\R^{d}} \left( u(x + A(C+x)^{1/2} yy^\top(C+x)^{1/2}A^\top) - u(x) \right)  \nu_L(dy) \\
&=: \mathcal D u(x) + \mathcal J u(x).
\end{align*}
We abbreviate the first addend, the drift part, with $\mathcal D u(x)$ and the second, the jump part, with $\mathcal J u(x)$. 

\paragraph{Foster-Lyapunov drift inequality}~\\
As test function we choose $u(x) = \tr(\eta x)^p + 1$, thus $u(x) \ge 1$. Note that the gradient of $u$ is given by $\nabla u(x) = p \tr(\eta x)^{p-1} \eta$. 

For $p \in(0,1)$ the gradient of $u$ has a singularity in $0$, but in the end we look at $\tr\left((Bx+xB^\top) \nabla u(x)\right)$, which turns out to be continuous in $0$. 

Now we need to look at the cases (i) - (v) separately, as the proofs work along similar lines, but differ in important details.
\subparagraph{Case (i)}~\\
We have
\begin{align*}
\mathcal J u(x) &= \int_{\R^{d}} \left( \tr\left(\eta(x + A(C+x)^{1/2} yy^\top(C+x)^{1/2}A^\top)\right) - \tr\left(\eta x\right) \right)  \nu_L(dy) \\
&=\int_{\R^{d}}  \tr\left( \eta^{1/2}A(C+x)^{1/2} yy^\top(C+x)^{1/2}A^\top\eta^{1/2}\right)  \nu_L(dy)\\
&=\tr\left( \eta^{1/2}A(C+x)^{1/2}\int_{\R^{d}}   yy^\top\nu_L(dy)(C+x)^{1/2}A^\top\eta^{1/2}\right)  \\
&\leq  \tr\left( \eta^{1/2}A(C+x)A^\top\eta^{1/2}\right)\left\|\int_{\R^{d}}  yy^\top  \nu_L(dy)\right\|_2.
\end{align*}
In the last step we used $X\preceq \|X\|_2I_d=\lambda_{max}(X)I_d$ for $X\in \s_d^+$, that $\tr$ is monotone in the natural order on $\s_d^+$,  and that maps of the form $\s_d\to\s_d, X\mapsto ZXZ^\top$ with $Z\in M_d(\R)$ are order preserving.

Hence
\begin{align*}
\mathcal J u(x) &\leq  \tr\left( A^\top\eta A(C+x)\right)\left\|\int_{\R^{d}}  yy^\top  \nu_L(dy)\right\|_2.
\end{align*}
For the drift part we get
\begin{align*}
\mathcal D u(x) & =\tr\left((Bx+xB^\top) \eta\right) =\tr((\eta B+B^\top \eta)x).
\end{align*}

Together we have
\begin{align*}
\mathcal A u(x) &= \tr\left(\left(\eta B+B^\top \eta+ A^\top\eta A\left\|\int_{\R^{d}}  yy^\top  \nu_L(dy)\right\|_2  \nu_L(dy)\right)x\right)+  \tr\left( A^\top\eta AC\right)\left\|\int_{\R^{d}}  yy^\top  \nu_L(dy)\right\|_2 .
\end{align*}
Since $ \eta B+B^\top \eta+ A^\top\eta A\left\|\int_{\R^{d}}  yy^\top  \nu_L(dy)\right\|_2\in-\s_d^{++}$, there exists a $c>0$, such that \[\tr\left(\left(\eta B+B^\top \eta+ A^\top\eta A\left\|\int_{\R^{d}}  yy^\top  \nu_L(dy)\right\|_2\right)x\right)\leq -c\tr(\eta x).\]

Summarizing we have that there exist $c,d>0$ such that
\begin{equation}
\label{gen1} \A u(x) \le - c \tr(\eta x) + d. 
\end{equation}

For $ \tr(x)> k$ and $k$ big enough there exist $0 < c_1 < c$ such that 
$
\A u(x)   \le - c_{1} ( \tr(\eta x) +1 ).
$
For $ \tr(x) \le k$ we have $
- c_1 \tr(\eta x) + d = -c_1 (\tr(\eta x) +1 ) + e,
$
with $e := c_1+d >0$.
Altogether we have
\begin{equation}  \label{gen2}
\A u(x) \le -c_1 u(x) + e \1_{D_{k}},
\end{equation}
where $D_{k}:= \{ x: \tr(x) \le k \}$ is a compact set. By Proposition \ref{prop:compactpetite} (i) this is also a petite set. Therefore the Foster-Lyapunov drift condition is proved. 

\subparagraph{Case (ii)}~\\
Compared to Case (i) we just change the inequality for  the jump part. Using Lemma \ref{lem:eigtrace} we have
\begin{align*}
\mathcal J u(x) 
&=\int_{\R^{d}}  \tr\left(A^\top\eta A(C+x)^{1/2} yy^\top(C+x)^{1/2}\right)  \nu_L(dy)\leq  \lambda_{max}(A^\top\eta A)   \tr\left((C+x) \int_{\R^{d}} yy^\top\nu_L(dy)\right) 
\end{align*}
and hence 
\begin{align*}
\mathcal A u(x) &= \tr\left(\left(\eta B+B^\top \eta+ \lambda_{max}(A^\top\eta A)\int_{\R^{d}} yy^\top\nu_L(dy)\right)x\right)+  \lambda_{max}(A^\top\eta A) \tr\left(C\int_{\R^{d}} yy^\top\nu_L(dy)\right).
\end{align*}
Now we can proceed as before.

\subparagraph{Cases (iii) and (iv)}~\\
We have for the drift part
\begin{align*}
\mathcal D u(x) & =p\tr((\eta B+B^\top \eta)x)\tr(\eta x)^{p-1}\leq K_{\eta,B} p \tr(\eta x)^p.
\end{align*}
For the jump part we get using that $\|yy^\top\|_2=\|y\|_2^2$ for $y\in \R^d$
\begin{align*}
\mathcal J u(x) &= \int_{\R^{d}} \left( \tr\left(\eta(x + A(C+x)^{1/2} yy^\top(C+x)^{1/2}A^\top)\right)^p - \tr\left(\eta x\right)^p \right)  \nu_L(dy) \\
&\leq\int_{\R^{d}} \left(\left( \tr\left(\eta x\right) + \|y\|^2_2\tr\left( \eta A(C+x)A^\top\right)\right)^p - \tr\left(\eta x\right)^p \right)\nu_L(dy)\\
&\leq\int_{\R^{d}} \left(\left( \tr\left(\eta x\right) + \|y\|^2_2\tr(\eta ACA^\top)+\|y\|^2_2K_{\eta,A}\tr(\eta x)\right)^p - \tr\left(\eta x\right)^p \right)\nu_L(dy).
\end{align*}
Using the elementary inequality $(x+y)^p \le \max\{2^{p-1},1\} (x^p + y^p)$ for all $x,y\geq 0$ we obtain
\begin{align*}\label{eq:jumpless1}
\mathcal J u(x) 
\leq&\tr\left(\eta x\right) ^p\int_{\R^{d}} \left(\max\{2^{p-1},1\} \left( 1+K_{\eta,A}\|y\|^2_2\right)^p - 1\right)\nu_L(dy)\\&+\max\{2^{p-1},1\}\tr(\eta ACA^\top)^p\int_{\R^{d}}  \|y\|^2_2\nu_L(dy).
\end{align*}
Putting everything together using \eqref{eq:assumptionp<1} or \eqref{eq:assumptionpge1}, respectively, we have that there exist $c,d>0$ such that
\begin{equation}
\A u(x) \le - c \tr(\eta x) ^p+ d. 
\end{equation}
and we can proceed as in Case (i).
\subparagraph{Case (v)}~\\
We again only argue differently in the jump part compared to the Cases (iii) and (iv) and can assume $p>1$ due to Remark \ref{rem44} (iii). For the jump part we again have
\[\mathcal J u(x) 
\leq\int_{\R^{d}} \left(\left( \tr\left(\eta x\right) + \|y\|^2_2\tr(\eta ACA^\top)+\|y\|^2_2K_{\eta,A}\tr(\eta x)\right)^p - \tr\left(\eta x\right)^p \right)\nu_L(dy).
\]
Next we use the following immediate consequence of the mean value theorem and the already used elementary inequality for $p$-th powers of sums.
\begin{lemma}
	For $x,y \ge 0$, $p \ge 1$ it holds that 
	\begin{equation*}
	0\leq (x+y)^p - x^p \le p y (x+y)^{p-1}.  
	\end{equation*}
\end{lemma}
This gives (using Landau $O(\cdot)$ notation)
\begin{align*}\mathcal J u(x) 
\leq&\int_{\R^{d}} \left(\left( \tr\left(\eta x\right) + \|y\|^2_2\tr(\eta ACA^\top)+\|y\|^2_2K_{\eta,A}\tr(\eta x)\right)^p - \tr\left(\eta x\right)^p \right)\nu_L(dy)\\
\leq& p \left(\tr(\eta ACA^\top)+K_{\eta,A}\tr(\eta x)\right)\\&\times
\int_{\R^{d}}\|y\|^2_2\left( \tr\left(\eta x\right) + \|y\|^2_2\tr(\eta ACA^\top)+\|y\|^2_2K_{\eta,A}\tr(\eta x)\right)^{p-1}\nu_L(dy)\\
=& pK_{\eta,A}\tr(\eta x)\int_{\R^{d}}\|y\|^2_2\left( \tr\left(\eta x\right) + \|y\|^2_2\tr(\eta ACA^\top)+\|y\|^2_2K_{\eta,A}\tr(\eta x)\right)^{p-1}\nu_L(dy)\\&+O(\tr(\eta x)^{p-1})\\
\leq& p\max\{2^{p-2},1\}K_{\eta,A}\tr(\eta x)^p\int_{\R^{d}}\|y\|^2_2\left(1 +\|y\|^2_2K_{\eta,A}\right)^{p-1}\nu_L(dy)+O\left(\tr(\eta x)^{\max\{p-1,1\}}\right).
\end{align*}
%
Combining this with the estimate on the drift part already given for Cases (iii) and (iv) and using  \eqref{eq:assumptionv}, we have that there exist $c,d>0$ such that
\begin{equation}
\A u(x) \le - c \tr(\eta x) ^p+ d. 
\end{equation}
and we can proceed as in Case (i).
\paragraph{Test function belongs to the domain}~\\
We now show that the chosen test function $u(x) = \tr(\eta x)^p +1$ belongs to the domain of the extended generator and that $\A u$ indeed has the claimed form. So we have to show that for all initial values $Y_0=x$ and all $t \ge 0$ it holds that
\begin{equation}
\label{exgen1_mucogarch} \E_{x}(u(Y_t)) = u(x) + \E_{x}\left( \int_{0}^{t}\A u(Y_s)ds\right) \\
\end{equation}
and
\begin{equation}
\label{exgen2_mucogarch} \E_{x}\left( \int_{0}^{t}|\A u(Y_s)|ds\right) < \infty.
\end{equation}
To show \eqref{exgen1_mucogarch} we show that 
$$ M_t := u(Y_t) - u(x) - \int_0^t \A u(Y_{s-}) ds$$
is a martingale. For this we apply It\^{o}'s formula for processes of finite variation to $u(x) = \tr(\eta x)^p +1$.  For a moment we ignore the discontinuity of $\nabla u$ in $0$ for $p \in (0,1)$ by considering $u$ on $\sdd\setminus\{0\}$. \\
\begin{align*}
u(Y_t) - u(Y_0) = &\int_{0+}^t \tr(\nabla u(Y_{s-}) d Y^c_s) + \sum_{0 < s \le t} \left(u(Y_s) - u (Y_{s-}) \right)\\
= & \int_{0+}^t p \zn{Y_{s-}}^{p-1} \tr\left((B Y_{s-}+Y_{s-}  B^\top) \eta\right) ds \\
&  + \int_{0+}^t \int_{\R^d} \left(u(Y_{s-} +  A(C+Y_{s-})^{1/2} yy^\top(C+Y_{s-})^{1/2}A^\top) - u(Y_{s-}) \right)\mu_{L}(dy,ds).
\end{align*}
Above we implicitly assumed that $\sum_{0 < s \le t} \left(u(Y_s) - u (Y_{s-}) \right)$ exists. Noting that $\tr(\eta Y_s )\geq \tr(\eta Y_{s-} )$ we easily see from the inequalities obtained for the jump part in the proof of the Foster-Lyapunov drift condition above  that this is the case whenever $L$ is of finite $2p$-variation which is ensured by the assumptions of the theorem.

Similarly we see that  \begin{equation*}\int_{0+}^t \int_{\R^d} \left(u(Y_{s-} +  A(C+Y_{s-})^{1/2} yy^\top(C+Y_{s-})^{1/2}A^\top) - u(Y_{s-}) \right)\nu_{L}(dy) ds\end{equation*}  is finite.
Then we get 
\begin{align*}
M_t &= u(Y_t) - u(x) - \int_{0}^{t}\A u(Y_{s-})ds \\
&= \int_{0+}^t \int_{\R^d} \left(u(Y_{s-} +  A(C+Y_{s-})^{1/2} yy^\top(C+Y_{s-})^{1/2}A^\top) - u(Y_{s-}) \right)(\mu_L(dy,ds)-\nu_L(dy)ds).
\end{align*} 
By the compensation formula (see \cite[Corollary 4.5]{kyprianou} for the version for conditional expectations), $(M_t)_{t \ge 0}$ is a martingale if \begin{equation}
\E \left[ \int_{0+}^t \int_{\R^d} \left(u(Y_{s-} +  A(C+Y_{s-})^{1/2} yy^\top(C+Y_{s-})^{1/2}A^\top) - u(Y_{s-}) \right)\nu_L (dy) ds \right] < \infty,
\label{eq:integrability}
\end{equation} as the integrand is non-negative.
Just like in the deduction of the Foster-Lyapunov drift condition, we get
$$\int_{\R^d}\left(u(Y_{s-} +  A(C+Y_{s-})^{1/2} yy^\top(C+Y_{s-})^{1/2}A^\top) - u(Y_{s-}) \right)\le c \zn{Y_{s-}}^p + d, $$
for some constants $c, d >0$. By Lemma \ref{lemma:moments} it holds that $\E(\zn{Y_t}^p) < \infty$ for all $t \ge 0$ and $ t \mapsto \E(\zn{Y_t}^p)$ is locally bounded. Thus \eqref{eq:integrability} follows by Fubini's theorem.

It remains to prove the validity of It\^{o}'s formula on the whole state space $\sdd$ and for $p\in(0,1)$. For this note that if $Y_0 \neq 0$ it follows that $Y_t \neq 0$ for all $t >0$. Thus we only need to consider the case $Y_0 = 0$. If the jumps are of compound Poisson type we define $\tau := \inf \{ t \ge 0: Y_t \neq 0 \}$, which is the first jump time of the driving L\'evy process. For $t < \tau$, $Y_t = 0$ and $t= \tau$ it is obviously fulfilled, because Lemma \ref{lem:etaB} implies that $C \tr(\eta Y_{s-})\geq|\tr\left((B Y_{s-}+Y_{s-}  B^\top) \eta\right) |\geq c \tr(\eta Y_{s-})$ for some $c,C>0$. For $t > \tau$ we can reduce it to the case $\sdd \setminus \{0\}$. Thus It\^{o}'s formula stays valid if the driving L\'evy process is a compound Poisson process. 

Now we assume that the driving L\'evy process has infinite activity. In  \cite[Proposition 6.7]{S2010} it has been shown, that we can approximate 
$Y_t$ by approximating the driving L\'evy process. We use the same notation as in \cite[Proposition 6.7]{S2010}. We fix $\omega \in \Omega$ and some $T >0$. The proof of \cite[Proposition 6.7]{S2010} shows that then $Y_{n,t}$ for all $n\in \mathbb{N}$ and $Y_{n,t}$ are uniformly bounded on $[0,T]$. It follows that $\zn{\cdot}^p$ is uniformly continuous on the space of values. We need to show that 
\begin{align*}
u(Y_{n,t}) - u(Y_0) 
= & \int_{0+}^t p \zn{Y_{n,s-}}^{p-1} \tr\left((B Y_{n,s-}+Y_{n,s-}  B^\top) \eta\right) ds \\
&  + \int_{0+}^t \int_{\R^d} \left(u(Y_{n,s-} +  A(C+Y_{n,s-})^{1/2} yy^\top(C+Y_{n,s-})^{1/2}A^\top) - u(Y_{n,s-}) \right)\mu_{L_n}(dy,ds).
\end{align*}
converges uniformly on $[0,T]$. Since $\omega \in \Omega$ and $T >0$ was arbitrary, once established this shows that It\^{o}'s formula holds almost surely uniformly on compacts. For the first summand uniform  convergence follows directly by the uniform convergence of $Y_{n,t}$ and the uniform continuity of the functions applied to it. 

 For the second summand we observe that also the jumps of $L$ are necessarily bounded on $[0,T]$ and that $$ \int_{0+}^t \int_{\R^d} \|y\|_2^{2p} \mu_{L_n}(dy,ds) \le \int_{0+}^t \int_{\R^d} \|y\|_2^{2p} \mu_{L}(dy,ds).$$ 
 We have that 
 \begin{align*}
 &\left|\int_{0+}^t \int_{\R^d} \left(u(Y_{n,s-} +  A(C+Y_{n,s-})^{1/2} yy^\top(C+Y_{n,s-})^{1/2}A^\top) - u(Y_{n,s-}) \right)\mu_{L_n}(dy,ds) \right.\\
  &\quad\quad\left.-\int_{0+}^t \int_{\R^d} \left(u(Y_{s-} +  A(C+Y_{s-})^{1/2} yy^\top(C+Y_{s-})^{1/2}A^\top) - u(Y_{s-}) \right)\mu_{L}(dy,ds) \right|\\
  &\leq \int_{0+}^T \int_{\R^d} \left|\left(u(Y_{n,s-} +  A(C+Y_{n,s-})^{1/2} yy^\top(C+Y_{n,s-})^{1/2}A^\top) - u(Y_{n,s-}) \right)\right.\\
  &\quad\quad\left.- \left(u(Y_{s-} +  A(C+Y_{s-})^{1/2} yy^\top(C+Y_{s-})^{1/2}A^\top) - u(Y_{s-}) \right) \right|\mu_{L}(dy,ds)\\
   &\quad\quad+\int_{0+}^T \int_{\R^d} \left|u(Y_{s-} +  A(C+Y_{s-})^{1/2} yy^\top(C+Y_{s-})^{1/2}A^\top) - u(Y_{s-}) \right|(\mu_{L}(dy,ds) -\mu_{L_n}(dy,ds)).
 \end{align*}
 Now the first integral converges to zero again by uniform continuity arguments for the integrand and the second one due to uniform boundedness of the integrand and the uniform convergence of $L_n$ to $L$.
 
It remains to show \eqref{exgen2_mucogarch}. 
To show this we first deduce a bound for $\abs{\A u(x)}$. Using the triangle inequality we split  $\abs{\A u(x)}$ again in the drift part and the jump part. For the absolute value of the jump part we can use the upper bounds from the Foster-Lyapunov drift condition proof since the jumps are non-negative. The absolute value of the drift part is bounded as follows
\begin{align*}
\abs{\mathcal D u(x)} &  \le pC \zn{Y_{n,s}}^p.
\end{align*}
Adding both parts together we get
\begin{equation} \label{gen3}
\abs{ \A(u(x)) } \le  c_2 u(x) 
\end{equation}
for some constant $c_2 >0$.

With that \eqref{exgen2_mucogarch} follows by Lemma \ref{lemma:moments}.\\
\paragraph{Harris recurrence and finiteness of moments}~\\
To show the positive Harris recurrence of the volatility process $Y$ and the finiteness of the $p$-moments of the stationary distribution we use the skeleton chains. In \cite[Theorem 5.1]{DMT1995} it is shown, that the Foster-Lyapunov condition for the extended generator, as we have shown, implies a Foster-Lyapunov drift condition for the skeleton chains. Further observe that petite sets are small, since by the assumption of irreducibility we can use the same arguments as in the upcoming proof of Theorem \ref{theo:irreducible}.  With that we can apply \cite[Theorem 3.12]{fuchs} and get the positive Harris recurrence for every skeleton chain and the finiteness of the $p$-moments of the stationary distribution. By definition the positive Harris recurrence for every skeleton chain implies it also for the volatility process $Y$.
\hfill $\Box$
 \subsection{Proof of Lemma \ref{lemma:moments}}$\E(\tr(Y_0)^{p}) < \infty$ implies $\E(\tr(\eta Y_0)^{p}) < \infty$ for all $\eta\in\s_d^{++}$.
 
 We apply It\^{o}'s formula to $u(Y_t) = \zn{Y_t}^p$. The validity of It\^o's formula was shown above in the section  ``Test function belongs to the domain''. We fix some $T >0$. Let $t \in [0,T]$. As in the proof before we get with It\^{o}'s formula 
 \begin{align*}
 u(Y_t) - u(Y_0) 
 = & \int_{0+}^t p \zn{Y_{s-}}^{p-1} \tr\left((B Y_{s-}+Y_{s-}  B^\top) \eta\right) ds \\
 &  + \int_{0+}^t \int_{\R^d} \left(u(Y_{s-} +  A(C+Y_{s-})^{1/2} yy^\top(C+Y_{s-})^{1/2}A^\top) - u(Y_{s-}) \right)\mu_{L}(dy,ds)
  \\
 \le & \int_{0+}^t c_1 \zn{Y_{s-}}^p ds  \\&+ \int_{0+}^t \int_{\R^d} \left(u(Y_{s-} +  A(C+Y_{s-})^{1/2} yy^\top(C+Y_{s-})^{1/2}A^\top) - u(Y_{s-}) \right)\mu_{L}(dy,ds)  .
 \end{align*}
 for some $c_1>0$.
 So we have 
 \begin{align*} 
 \E(\zn{Y_t}^p) 
 \le & \E(\zn{Y_0}^p) + \int_0^t c_1 \E(\zn{Y_s}^p) ds \\&+ \E \left( \int_{0+}^t \int_{\R^d} \left(u(Y_{s-} +  A(C+Y_{s-})^{1/2} yy^\top(C+Y_{s-})^{1/2}A^\top) - u(Y_{s-}) \right)\mu_{L}(dy,ds) \right).
 \end{align*}
 Using the compensation formula and the bounds of the proof of the Foster-Lyapunov drift condition we get that there exist $c_2,d>0$ such that 
 \begin{align*}&\E \left( \int_{0+}^t \int_{\R^d} \left(u(Y_{s-} +  A(C+Y_{s-})^{1/2} yy^\top(C+Y_{s-})^{1/2}A^\top) - u(Y_{s-}) \right)\mu_{L}(dy,ds)\right)\\&\quad\quad
 \le \int_0^t \E( c_2 \zn{Y_s}^p + d) ds . 
 \end{align*}
 Combined we have 
 \begin{align*}
 \E(\zn{Y_t}^p) & \le \E(\zn{Y_0}^p) + \int_0^t (c_1+c_2) \E(\zn{Y_s}^p) ds  + d T. 
 \end{align*}
 Applying Gronwall's inequality this shows that
 $$ \E(\zn{Y_t}^p) \le ( \E(\zn{Y_0}^p) + d T) e^{(c_1+c_2) t}. $$
 Since $T >0$ was arbitrary  $\E(\zn{Y_t}^p)$ is finite for all $t \ge 0$ and $t \mapsto \E(\zn{Y_t}^p)$ is locally bounded.

\subsection{Proofs of Section \ref{section:irredresults}} 
\subsubsection{Proof of Theorem \ref{theo:irreducible}} \label{sectionirreducibility}

Let $\nu$ be the L\'evy measure of $L$. We have $\nu= \nu_{\text{ac}} + \tilde \nu$, where $\nu_{\text{ac}}$ is the absolute continuous component and $\nu(\R) < \infty$. Moreover we can split $L$ into the corresponding processes $L\stackrel{\mathcal D}{=} L_{\text{ac}}+\tilde L$, where $L_{\text{ac}}$ and $\tilde L$ are independent.  We set $B_T=\{\omega \in \Omega~ | ~ \tilde L_t = 0~ \forall t \in [0,T]\}$. Then $\forall T >0~~ \p(B_T)>0$ and for any event $ A$ it holds that
$\p(A \cap B_T) >0 \Leftrightarrow \p(A|B_T) >0$
and $\p(A\cap B_T) >0 \Rightarrow \p(A) >0$. 
So in the following we assume w.l.o.g. $\tilde L = 0$ as otherwise the below arguments and the independence of $L_{\text{ac}}$ and $\tilde L$ imply that $\p(Y_t \in A~|~Y_0=x)>0$ results from $\p(Y_t \in A ~|~ Y_0 =x, B_t) >0$.

\paragraph*{Irreducibility:}~ \\
By Remark \ref{irredskeleton}, to prove the irreducibility of 
the MUCOGARCH volatility process it is enough to 
show it for a skeleton chain.  

Let $\delta >0$ and set $t_n := n\delta, \forall n\in \N_0$.  We consider the skeleton chain 
\begin{align*}
 Y_{t_n} =& e^{Bt_n } Y_{t_0} e^{B^\top t_n} + \int_0^{t_n} e^{B(t_n-s)} A 
(C+Y_{s-})^{\frac 1 2} d[L,L]^d_s (C+Y_{s-})^{\frac 1 2} A^\top e^{B^\top 
(t_n-s)}.  \\
\end{align*}

To show irreducibility w.r.t. $\lambda_{\sdd}$ we have to show that for any $\mathcal A \in \B(\s_{d}^{+})$ with $\lambda_{\s_d^+}(\mathcal A) >0$ and any $y_0 \in \sdd$ there exists an $l$ such that 
\begin{equation}
 \p(Y_{t_{l}} \in \mathcal A | Y_{0}=y_0) >0. \qquad 
\end{equation}
With
\begin{align}
\nonumber &\p(Y_{t_{l}} \in \mathcal A | Y_{t_0}=y_0) \\
 \nonumber &\ge \p(Y_{t_{l}} \in \mathcal A, \text{ exactly one jump in every 
time interval }(t_0,t_{1}], \cdots, (t_{l-1},t_{l}] | Y_{t_0}=y_0) \\
\nonumber &=  \p(Y_{t_{l}} \in \mathcal A | Y_{t_0}=y_0,\text{ exactly one 
jump in every time interval }(t_0,t_{1}], \cdots, (t_{l-1},t_{l}] ) \\
\label{wkeitonejump}&\cdot\p(\text{ exactly one jump in every time interval 
}(t_0,t_{1}], \cdots, (t_{l-1},t_{l}]), 
\end{align}
and the fact that the last factor is strictly positive we can w.l.o.g. assume, that we have exactly one jump in every time interval $(t_{k-1}, t_{k}] \ \forall k=1, \cdots, l$.

We denote by $\tau_k$ the jump time of our L\'evy process in $(t_{k-1}, t_{k}]$. With the assumption, that we only have one jump on every time interval, the skeleton chain can be represented by the sum of the jumps $X_i$, where $L_t=\sum_{i=1}^{N_t} X_i$ is the used representation for the compound Poisson Process $L$. We fix the number of time steps $l\ge d$ and get:
\begin{align}
 \label{ylstep} Y_{t_{l}}=~& e^{Bt_{l}} Y_{0} e^{B^\top t_{l} }  \\&+ \sum_{i=1}^{l} e^{B(t_{l}- \tau_i)} A (C+e^{B(\tau_i-t_{i-1})}Y_{t_{i-1}}e^{B^\top(\tau_i-t_{i-1})})^{\frac 1 
2} X_{i}X_{i}^\top (C+e^{B(\tau_i-t_{i-1})}Y_{t_{i-1}}e^{B^\top(\tau_i-t_{i-1})})^{\frac 1 2} A^\top 
e^{B^\top(t_{l}-\tau_i)}.\nonumber
\end{align}

First we show that the sum of jumps in (\ref{ylstep}) has a positive density on $\s_d^+$. 
Note that every single jump is of rank one. 

We define 
\begin{equation}
 \label{zi} Z_i^{(l)} := e^{B(t_{l}- \tau_i)} A (C+e^{B(\tau_i-t_{i-1})}Y_{t_{i-1}}e^{B^\top(\tau_i-t_{i-1})})^{\frac 1 2} X_{i}
\end{equation}
  
and with (\ref{ylstep}) we have 

\begin{align*}
 Z_i^{(l)} = e^{B(t_{l}- \tau_i)} A \Big( C &+ e^{Bt_{i-1}} Y_{0} e^{B^\top t_{i-1}} + \sum_{j=1}^{i-1} e^{B(t_{i-1}-t_{l})} Z_j^{(l)} {Z_j^{(l)}}^\top e^{B^\top (t_{i-1}-t_{l})} \Big)^{\frac 1 2} X_{i}.
\end{align*} 

By assumption $X_1, X_2, \cdots$ are  are iid and absolutely continuous w.r.t. Lebesgue measure on $\R^d$ with a (Lebesgue a.e.) strictly positive density. We see 
immediately that 
$
 Z_1^{(l)} | Y_{0}, \tau_1
$
is absolutely continuous with a strictly positive density $f_{Z_1^{(l)} | Y_{0}, \tau_1}$. Iteratively we get that 
every
\begin{equation}
 Z_i^{(l)} | Y_{0}, Z_1^{(l)}, \cdots, Z_{i-1}^{(l)}, \tau_i
\end{equation}
is absolutely continuous with a strictly positive density $f_{Z_i^{(l)} | Y_{0}, Z_1^{(l)}, \ldots, Z_{i-1}^{(l)}, \tau_i }$, 
for all $i=2, \ldots, l$.  

We denote with $f_{Z^{(l)}|Y_{0}, \tau_1, \ldots, \tau_l}$ the density of $Z^{(l)}=(Z_1^{(l)}, \cdots, Z_l^{(l)})^\top$ 
given $Y_{0}, \tau_1, \ldots, \tau_l$. Note that given $Z_j^{(l)},~j<i$, $Z_i^{(l)}$ is independent of $\tau_j$. By the rules for conditional densities we get
\begin{equation*}
 f_{Z^{(l)}|Y_{0}, \tau_1, \ldots, \tau_l} = f_{Z_1^{(l)}| Y_{0},\tau_1} \cdot f_{Z_2^{(l)} | Y_{0}, Z_1^{(l)}, \tau_2} \cdots f_{Z_l^{(l)}| Y_{0}, Z_1^{(l)}, \cdots, Z_{l-1}^{(l)}, \tau_l}
\end{equation*}
is strictly positive on $\R^{dl}$. Thus an equivalent measure $\q$, $\q \sim 
\p$, exists such that \linebreak $Z_1^{(l)}|Y_{0},\tau_1,\ldots, \tau_l, ~ \cdots ~, Z_l^{(l)}|Y_{0}, \tau_1,\ldots. \tau_l$ are iid normally distributed. In \cite{anderson} it is shown, that for $l \ge d$ 
\begin{equation}
\Gamma:= \sum_{i=1}^l Z_i^{(l)}|_{Y_{0}, \tau_1,\ldots. \tau_l} \cdot {Z_i^{(l)}}^\top |_{Y_{0}, \tau_1,\ldots. \tau_l}
\end{equation}
has a strictly positive density under $\q$ w.r.t. the Lebesgue measure 
on $\s_d^+$. But since $\q$ and $\p$ are equivalent, $\Gamma$ has also a 
strictly positive density under $\p$ w.r.t. Lebesgue measure 
on $\s_d^+$. 

This yields 
\begin{align*}
&\p(Y_{t_l} \in \mathcal A | Y_0=y_0) =\int_{\R_+^{l}} \p( e^{Bt_l}Y_0e^{B^\top t_l} + \Gamma \in \mathcal A|Y_0=y_0, \tau_1=k_1, \ldots, \tau_l=k_l) d\p_{(\tau_1,\ldots, \tau_l)}(k_1,\ldots, k_l) >0
\end{align*}
if $\p( e^{Bt_l}Y_0e^{B^\top t_l} + \Gamma \in \mathcal A~|~Y_0=y_0, \tau_1=k_1, \ldots, \tau_l=k_l)>0.$ Here we use that the joint distribution of the jump times $\p_{(\tau_1,\ldots, \tau_l)}$ is not trivial, since $\tau_1, \ldots, \tau_k$ are the jump times of a compound Poisson Process. 
Above we have shown that $\p( e^{Bt_l}Y_0e^{B^\top t_l} + \Gamma \in \mathcal A~|~Y_0=y_0, \tau_1=k_1, \ldots, \tau_l=k_l)>0$ if \[\lambda_{\s_d^+} \left(\mathcal A \cap \{x \in \s_d^+ | x \succeq e^{B t_{l}}Y_{0} e^{B^\top t_{l}} \} \right) >0.\] 

As we assumed $\sigma(B) \subset (-\infty,0) +i\R$,  $e^{Bt}Y_{0}e^{B^\top t} \to 0$  for $t \to \infty$. Thus we can choose $l$ big enough such that $\lambda_{\s_d^+} \left(\mathcal A \cap \{x \in \s_d^+ | x \succeq e^{B t_{l}}Y_{0} e^{B^\top t_{l}} \} \right) >0$
for any $\mathcal A \in \s_d^+$ with $\lambda_{\s_d^+}(\mathcal A) >0$ (note  $\lambda_{\s_d^+}(\partial \s_d^+) =0$). This shows the claimed irreducibility  and as $\delta$ was arbitrary even simultaneous irreducibility. 

\paragraph*{Aperiodicity:}~ \\
The 
simultaneous irreducibility and Proposition 
\ref{prop:simirred} show that every skeleton chain is aperiodic. Using Proposition \ref{prop:compactpetite} we know for every skeleton chain, that every compact set is also small.

We define the set \begin{equation}
                     \mathcal{C} := \{ x \in \s^d_+ | \ \|x\|_2 \le K\},
                    \end{equation}
with a constant $K>0$. Obviously $\mathcal{C}$ is a compact set and thus a small set for every skeleton chain. By Remark \ref{remark:small} it is also small for the continuous time Markov Process $(Y_t)_{t\ge 0}$.
To show aperiodicity for $(Y_t)_{t\ge 0}$ in the sense of Definition \ref{def:aperiodic} we prove that 
there exists a $T>0$ such that
\begin{equation} \label{aperiodictoshow}
 \p^t(x,\mathcal{C}) >0
\end{equation}
holds for all $x\in \mathcal C$ and all $t\ge T$. 

Using 
\begin{equation} \label{eq:aperiodicnojump}
 \p^t(x,\mathcal C) \ge \p^t(x, \mathcal C \cap \{\text{ no jump up to time 
$t$}\})
\end{equation}
we consider $Y_t$ under the condition ``no jump up to time $t$''. With
$Y_0 =x \in \mathcal C$ we have
\begin{equation} \label{eq:Ynojump}
 Y_t = e^{Bt}xe^{B^\top t}.
\end{equation}
Since $\lambda = \max ( \re(\sigma(B))) <0$ there exists $\delta >0$,  and $C 
\ge 1$ such that $ \| e^{Bt} \|_2 \le C e^{- \delta t}$ and hence we have
$
       \| e^{Bt} x e^{B^\top t} \|_2 \le C e^{-2\delta t} \|x\|_2\le C e^{-2\delta t} K  \le K
$
for all $t \ge  \frac{\ln(C)}{2\delta}$. Hence, $Y_t \in \mathcal C$ for all $t \ge  \frac{\ln(C)}{2\delta}$ and thus (\ref{aperiodictoshow}) holds. \hfill $\Box$

\subsubsection{Proof of Corollary \ref{cor:irredcppnullenvironment}}
The proof is similar to that of Theorem \ref{theo:irreducible} with the difference that we now assume for the jump sizes $X_i$ that they have a density, which is strictly positive in a neighborhood of zero, e.g. $\exists k>0$ such that 
every $X_i$ has a strictly positive density on $\{ x \in \R^d: \| x \| \le k\}$. 
We use the same notation as in the previous proof, but we omit the superscripts $^{(l)}$. 
By the definition of $Z_i$ and the same iteration as in the first case we show, that $Z:=(\tilde Z_1, \ldots, \tilde Z_l)|Y_{0},\tau_1, \ldots, \tau_l$ has a strictly positive density on a suitable neighborhood of the origin.

For $A\in M_d(\R)$ we define $j(A):= \min_{x \in \R^d} \frac{\|Ax\|_2}{\|x\|_2}$ as the modulus of injectivity, which has the following properties: $0 \le j(A) \le \| A\|_2$ and $ \| Ax\|_2 \ge j(A) \|x\|_2$ as well as for $A,B \in M_d(\R)$  $j(A  B) \ge j(A)j(B)$. \\
With that we get for $Z_1$
\begin{align*}
\| Z_1 \|_2 &= \| e^{B(t_l-\tau_1)} A (C+e^{B\tau_1}Y_0e^{B^\top\tau_1})^{\frac 1 2 } X_1\|_2  \ge j(e^{B(t_l-\tau_1)} A (C+e^{B\tau_1}Y_0e^{B^\top\tau_1})^{\frac 1 2 }) ~ \|X_1\|_2 \\
& \ge j(e^{B(t_l-\tau_1)})~ j( A)~ j( (C+e^{B\tau_1}Y_0e^{B^\top\tau_1})^{\frac 1 2 }) ~\|X_1\|_2  \ge j(e^{Bt_l})~ j( A)~ j( C^{\frac 1 2 }) ~\|X_1\|_2 
\end{align*}
and thus $Z_1|Y_{0},\tau_1$ has a strictly positive density on $\{x\in \R^d: \|x\| \le \tilde k\}$, where $\tilde k := j(e^{Bt_l})~ j( A)~ j( C^{\frac 1 2 }) ~ k$. 
Iteratively get that every $Z_i|Y_0, \tau_i, Z_1, \ldots, Z_{i-1}$ has a strictly positive density on $\{x\in \R^d: \|x\| \le \tilde k\}$ and as in the first case this shows that $Z=( Z_1, \ldots,  Z_l)|Y_{0},\tau_1, \allowbreak \ldots,\allowbreak \tau_l$ has a strictly positive density on $\{ x=(x_1,\ldots,x_l)^\top \in \R^{d \cdot l}: \| x_i \| \le \tilde k~ \forall i=1,\ldots,l \}.$

We fix an $\hat k$, $0<\hat k < \tilde k$ and set $\mathcal{\hat K}:= \{x =(x_1,\ldots,x_l)^\top\in \R^{d \cdot l}:\|x_i\| \le \hat k \  \forall i=1, \ldots,l \}$. Now we can construct random variables $\tilde Z_i$, $i=1, \ldots, l$, such that 
 $\tilde Z:=(\tilde Z_1, \ldots, \tilde Z_l)|Y_{0},\tau_1, \ldots, \tau_l$ has a strictly 
positive 
density on $\R^{d \cdot l}$ 
and 
\begin{equation}
\1_{\mathcal{\hat K}} \cdot \tilde Z |Y_{0},\tau_1, \ldots, \tau_l ~  \stackrel{\mathcal D}{=} \1_{\mathcal{\hat K}} \cdot  Z |Y_{0},\tau_1, \ldots, \tau_l .
\end{equation}

Due to the first case we now can choose a measure $\q$ such that the $\tilde 
Z_i|_{Y_0, \tau_1, \ldots, \tau_l, Z_1, \ldots, Z_{i-1}}$ are iid normal distributed and the random variable $\tilde \Gamma := \sum_{i=1}^l \tilde Z_i|_{Y_0, \tau_1, \ldots, \tau_l, Z_1, \ldots, Z_{i-1}} ~  \tilde Z_i^\top |_{Y_0, \tau_1, \ldots, \tau_l, Z_1, \ldots, Z_{i-1}}$ has a strictly positive density on $\s_d^+$. \\
With the equivalence of $\p$ and $\q$ also $\p(\tilde \Gamma \in \mathcal A)>0$ for every $\mathcal A \in \mathcal B(\sdd)$.

Further we define $\mathcal E:=\{x\in \sdd: x=\sum_{i=1}^l z_i z_i^\top,\ z_i\in \R^d, \|x\|_2 \le \hat k\}$ and $\mathcal K:= \{x\in \sdd :   x=\sum_{i=1}^l z_i z_i^\top \text{ for }\allowbreak z_1, \ldots, z_l \in \R^d  \text{ implies } \|z_i\|_2 \le \hat k\ \forall i=1,\ldots,l \}$. Let $x= \sum_{i=1}^l z_i z_i^\top \in \mathcal E$. Then $x = \sum_{i=1}^l z_i z_i^\top \succeq z_j z_j^\top$ for all $j=1,\ldots, l$ and thus $\|z_j z_j^\top\|_2=\|z_j\|_2\le \hat k$, which means that $x \in \mathcal K$ and thereby $\mathcal E \subseteq \mathcal K$. 

Now let $\mathcal A \in \mathcal B (\sdd)$ and note that $\1_{\mathcal{K}} \cdot \tilde \Gamma ~  \stackrel{\mathcal D}{=} \1_{\mathcal{ K}} \cdot  \Gamma$. Finally we get 
\begin{align*}
 \p(\Gamma \in \mathcal A \cap \mathcal E) &= \p(\Gamma \in \mathcal{A} \cap 
\mathcal{E} \cap \mathcal{K}) = \p(\tilde \Gamma \in \mathcal{A} \cap \mathcal{E} )>0
\end{align*}
if $\lambda_{\s_d^+}(\mathcal{A} \cap \mathcal{E})>0$.  

With the same conditioning argument as in the proof of Theorem \ref{theo:irreducible} and again using the fact that by assumption there always exists a $l$ such that $\lambda_{\s_d^+}(\mathcal{A} \cap \mathcal{E} \cap \{x \in \s_d^+ | x \succeq e^{B t_{l}}Y_{0} e^{B^\top t_{l}} \})>0$ if $\lambda_{\s_d^+}(\mathcal{A} \cap \mathcal{E})>0$, we get simultaneous irreducibility w.r.t. the measure $\lambda_{\s_d^+ \cap \mathcal E}$ defined by $\lambda_{\s_d^+ \cap \mathcal E}(B) := \lambda_{\s_d^+}(B \cap \mathcal{E})$ for all $B \in \mathcal B (\sdd)$. 

Aperiodicity follows as in the proof of Theorem \ref{theo:irreducible}, since we only used the compound Poisson structure of $L$ and not the assumption on the jump distribution. \hfill $\Box$

\section*{Acknowledgements}
The authors are grateful to the editor and an anonymous referee for their reading and insightful comments which improved the paper considerably.
The second author gratefully acknowledges support by the DFG Graduiertenkolleg 1100.

\bibliographystyle{acmtrans-ims1}
{\small
	\providecommand{\MR}[1]{}

	\bibliography{Literatur}}

\end{document}